\documentclass[preprint,sort&compress,final,10pt]{elsarticle}
\usepackage{amsmath}
\usepackage{amssymb}
\usepackage{graphicx}
\usepackage{latexsym}
\usepackage{epstopdf}
\usepackage{natbib}
\usepackage{color}
\usepackage{hyperref}
\usepackage{amsmath,amscd}
 \usepackage[all]{xy}

\newtheorem{theorem}{Theorem}[section]

\newtheorem{definition}[theorem]{Definition}
\newtheorem{lemma}[theorem]{Lemma}

\newtheorem{remark}[theorem]{Remark}
\numberwithin{equation}{section}

\def\Proof{\noindent{\bf Proof.}~}
\def\qed{\hfill$\square$\smallskip}

\def\dsum{\displaystyle\sum}

\def\Re{\mathrm{Re}}
\def\Im{\mathrm{Im}}

\journal{\empty}
\date{}

\begin{document}

\begin{frontmatter}

\title{Invariant curves of smooth quasi-periodic  mappings}

\author[au1]{Peng Huang}

\address[au1]{School of Mathematics Sciences, Beijing Normal University, Beijing 100875, P.R. China.}

\ead[au1]{hp@mail.bnu.edu.cn}

\author[au1]{Xiong Li\footnote{Corresponding author. Partially supported by the NSFC (11571041) and the Fundamental Research Funds for the Central Universities.}}

\ead[au1]{xli@bnu.edu.cn}

\author[au2]{Bin Liu\footnote{Partially supported by the NSFC (11231001).}}

\address[au2]{School of Mathematical Sciences, Peking University, Beijing 100871, P.R. China.}

\ead[au3]{bliu@pku.edu.cn}

\begin{abstract}
In this paper we are concerned with the existence of invariant curves of planar mappings which are quasi-periodic in the spatial variable, satisfy the intersection property, $\mathcal{C}^{p}$ smooth with $p>2n+1$, $n$ is the number of frequencies.
\end{abstract}


\end{frontmatter}

\section{Introduction}
In this paper we are  concerned with the  existence of  invariant curves of the following planar quasi-periodic  mappings
\begin{equation}\label{M}
\mathfrak{M}:\quad \begin{array}{ll}
\left\{\begin{array}{ll}
\theta_1=\theta+r +f(\theta,r),\\[0.2cm]
r_1=r+g(\theta,r),\\[0.1cm]
 \end{array}\right.\  (\theta,r)\in \mathbb{R} \times [a,b],
\end{array}
\end{equation}
where the perturbations $f(\theta,r)$ and $g(\theta,r)$ are quasi-periodic in $\theta$ with the frequency $\omega=(\omega_1,\omega_2$, $\cdots,\omega_n)$,\ $\mathcal{C}^p$ smooth  in $\theta$ and $r$.\

In 1962,\ Moser \cite{Moser62} considered  the twist mapping
$$
\mathfrak{M}_{0}:\quad \begin{array}{ll}
\left\{\begin{array}{ll}
x_1=x+\alpha(y)+\varphi_{1}(x,y),\\[0.2cm]
y_1=y+\varphi_{2}(x,y),
 \end{array}\right.
\end{array}
$$
where the perturbations $\varphi_{1},\varphi_{2}$ are assumed to be small and of periodic $2\pi$ in $x.$\ He obtained the existence of invariant closed curves of  $\mathfrak{M}_{0}$ which is of class $\mathcal{C}^{333}$.\ About  $\mathfrak{M}_{0}$,\ an analytic version of the invariant curve theorem was presented in \cite{Siegel97},\ a version in class $\mathcal{C}^{5}$ in  R\"{u}ssmann \cite{Russmann70} and a optimal version in class $\mathcal{C}^{p}$ with $ p>3$  in Herman \cite{Herman83, Herman86}.

When the perturbations $f(\theta,r),g(\theta,r)$ in (\ref{M}) are quasi-periodic in $\theta$, there are some results about the existence of invariant curves of the following planar quasi-periodic  mappings
\begin{equation}\label{a7}
\mathfrak{M}_{1}:\quad \begin{array}{ll}
\left\{\begin{array}{ll}
\theta_1=\theta+\beta+r +f(\theta,r),\\[0.2cm]
r_1=r+g(\theta,r),
 \end{array}\right.\ \ \ \   (\theta,r)\in \mathbb{R} \times [a,b],
\end{array}
\end{equation}
where the functions $f(\theta,r)$ and $g(\theta,r)$ are quasi-periodic in $\theta$ with the frequency $\omega=(\omega_1,\omega_2$,$\cdots,\omega_n)$, real analytic in  $\theta$ and $r$, and $\beta$ is a constant.

When the map $\mathfrak{M}_{1}$ in (\ref{a7}) is an exact symplectic map,\ $\omega_1,\omega_2,\cdots,\omega_n$, $2\pi\beta^{-1}$ are sufficiently incommensurable, Zharnitsky \cite{Zharnitsky00} proved the existence of invariant curves of the map $\mathfrak{M}_{1}$ and applied this result to present the boundedness of all solutions of Fermi-Ulam problem.\ His proof is based on the Lagrangian approach introduced by Moser \cite{Moser88} and used by Levi and Moser in \cite{Levi01} to show a proof of the twist theorem.

When the map $\mathfrak{M}_{1}$ in (\ref{a7}) is  reversible with respect to the involution $\mathcal{G}:(x,y)\mapsto (-x,y)$,\ that is,\ $\mathcal{G} \mathfrak{M}_1 \mathcal{G} = \mathfrak{M}_1^{-1}, \omega_1,\omega_2,\cdots,\omega_n,2\pi\beta^{-1}$ satisfy the Diophantine condition
\begin{equation*}
\Big|{\langle k,\omega \rangle {\beta \over {2\pi}}-j}\Big|\geq {\gamma \over {|k|^\tau}},\ \ \ \ \forall\ \  k \in \mathbb{Z}^n\backslash\{0\},\ \ \forall j \in \mathbb{Z},
\end{equation*}
Liu  \cite{Liu05} stated some variants of the invariant curve theorem for quasi-periodic reversible mapping $\mathfrak{M}_{1}$.

In this paper,\ motivated by the above references,\ especially by R\"{u}ssmann \cite{Russmann83}, instead of the exact symplecticity or reversibility assumption on $\mathfrak{M}$,\ we assume that this mapping satisfies the intersection property,\ and  obtain the  invariant curve theorem for  the quasi-periodic  mapping $\mathfrak{M}$ in the smooth case, other than analytic case.  

Incidently, in \cite{Huang16} we use this theorem to establish the  existence of invariant curves of the  planar quasi-periodic  mapping 
\begin{equation*}
\begin{array}{ll}
\mathcal{M}_{\delta}:\ \ \left\{\begin{array}{ll}
\theta_1=\theta+\beta+\delta l(\theta,r)+\delta f(\theta,r,\delta),\\[0.2cm]
r_1=r+\delta m(\theta,r)+\delta g(\theta,r,\delta),
 \end{array}\right.\ \ \  (\theta,r)\in \mathbb{R} \times [a,b],
\end{array}
\end{equation*}
where the functions $l,m, f, g$ are quasi-periodic in $\theta$ with the frequency $\omega$=$(\omega_1$,$\omega_2$,\\ $\cdots$, $\omega_n),$
 $f(\theta,r,0 ) =g(\theta,r, 0 ) = 0, $\ $\beta$ is a constant, $0<\delta< 1$ is a small parameter. As an application, we also use them to study the existence of quasi-periodic solutions and the boundedness of all solutions for an asymmetric oscillation 
\begin{equation*}
\begin{array}{ll}
x''$+$a$$x^+$$- b$$x^-$=$f(t),
\end{array}
\end{equation*}
where $a$, $b$ are two different positive constants,\ $x^+=\max\{x,0\}$,\ $x^-=\max\{-x,0\}$,\ $f(t)$ is  a smooth quasi-periodic function with the frequency $\omega=(\omega_1,\omega_2,\cdots, \omega_n)$.

Finally, we must point out that in order to obtain the existence of invariant curves for the quasi-periodic  mapping $\mathfrak{M}$, we need to assume that this mapping belongs to $\mathcal{C}^{p}$ with $p>2n+1$ and $n$ is the number of the frequency $\omega=(\omega_1,\omega_2,\cdots,\omega_n)$. Meanwhile we note that when $n=1$, quasi-periodic mappings are periodic mappings, and the optimal smoothness assumption is $\mathcal{C}^{p}$ with $p>3$. Hence our smoothness assumption for quasi-periodic mappings agrees with that for periodic mappings, and is optimal in this sense.

Our efforts in this paper are same as R\"{u}ssmann \cite{Russmann83}, we are more interested in weak conditions for the perturbations $f,g$ than in high differentiability properties of the constructed invariant curves, and the main line of the proofs is also similar to that of R\"{u}ssmann \cite{Russmann83}.

The rest of the paper is organized as follows.\ In Section 2,\ we  list some properties of quasi-periodic functions,\ and then  state the main invariant curve theorem  (Theorem \ref{thm2.6}) for  the quasi-periodic mapping $\mathfrak{M}$ which is given by (\ref{M}).\ The proofs of Theorem \ref{thm2.6} are given in Sections 3,\ 4,\ 5. In this section 6, we formulate the detail proofs of the Lemma \ref{lem2.9} which have been used in the previous sections.

\section{Quasi-periodic functions and the main result}
\subsection{The space of  quasi-periodic functions}

We first define quasi-periodic functions with the frequency $\omega$, here the $n$-dimensional frequency vector $\omega=(\omega_1,\omega_2,\cdots,\omega_n)$ is rationally independent, that is, for all $k=(k_1, k_2,\cdots, k_n) \neq 0$,\ $\langle k,\omega \rangle =\dsum \limits_{j=1}^n k_j \omega_j \neq 0$.

\begin{definition}\label{def2.1}
$f(t)$ is called a continuous quasi-periodic function with the frequency $\omega=(\omega_1,\omega_2,\cdots,\omega_n)$, if there is a continuous function $F(\theta_1,\theta_2,\cdots,\theta_n)$ which is $2\pi$-periodic in each $\theta_j\, (1\leq j\leq n)$ such that
$$f(t)=F(\omega_1 t,\omega_2 t,\cdots,\omega_n t).$$
Moreover, $f(t)$ is called a $\mathcal{C}^p $/real analytic quasi-periodic function, if $F$ is $\mathcal{C}^p $/real analytic, meanwhile we say that $F$ is a shell function of $f$.
\end{definition}

Denote by $Q(\omega)$ the space of real analytic quasi-periodic functions with the frequency $\omega=(\omega_1,\omega_2,\cdots,\omega_n)$. Given
$f(t)\in Q(\omega)$, suppose that the corresponding shell  function $F$ has the following Fourier expansion
$$F(\theta)=\sum _{k\in \mathbb{Z}^n} f_{k}e^{i \langle k,\theta \rangle },$$
which is $2\pi$-periodic in each variable, real analytic and bounded in a complex neighborhood $\Pi_{r}^{n}=\{(\theta_{1},\theta_{2},\cdots,\theta_{n})\in \mathbb{C}^n : |\Im\ \theta_{j}|\leq r, j=1,2,\cdots, n \}$  of $\mathbb{R}^n$ for some $r > 0$.\ The function $f(t)$ is obtained from $F(\theta)$ by replacing $\theta$ by $\omega t$, and has the following expansion
$$f(t)=\sum \limits _{k\in \mathbb{Z}^n} f_{k}e^{i \langle k,\omega \rangle t}.$$

\begin{definition}\label{def2.2}
For $r>0$, let $Q_{r}(\omega)\subseteq Q(\omega)$ be the set of real analytic quasi-periodic functions $f$ such that the corresponding shell functions $F$ are bounded on the subset $\Pi_{r}^{n}$\ with the supremum norm
$$\big|F\big|_{r}=\sup \limits_{\theta\in \Pi_{r}^{n}}|F(\theta)|=\sup \limits_{\theta\in \Pi_{r}^{n}}\Big|\sum_{k\in \mathbb{Z}^n}f_{k}e^{i\langle k,\theta\rangle}\Big|<+\infty.$$
Also we define the norm of $f$ as $\big|f\big|_{r}=\big|F\big|_{r}.$
\end{definition}

The following properties of quasi-periodic functions can be found in \cite[chapter 3]{Siegel97}.

\begin{lemma}\label{lem2.3}
The following statements are true: \\[0.2cm]
$(i)$ Let $f(t),g(t)\in Q(\omega)$,\ then $g(t+f(t))\in Q(\omega);$\\[0.2cm]
$(ii)$ Suppose that
\begin{equation*}
\begin{array}{ll}
{{| \langle k,\omega \rangle|} \geq {c \over {|k|^{\sigma_{0}}}}},\ \ \ \ c,\sigma_{0} > 0
\end{array}
\end{equation*}
for all integer vectors $k\neq0$.\ Let $h(t)\in Q(\omega)$ and $\tau=\beta t +h(t)\ (\beta+h^{'}>0)$,\\[0.2cm] then the inverse relation is given by $t={\beta^{-1}} \tau +h_1(\tau)$  and $h_1\in Q({\omega \over \beta})$.\ In particular,\ if $\beta=1$,\ then $h_1\in Q(\omega).$
\end{lemma}

Throughout this paper,\ we assume that the frequency  $\omega=(\omega_1,\omega_2,\cdots,\omega_n)$ satisfies the Diophantine condition
\begin{equation}\label{b2}
\begin{array}{ll}
{{| \langle k,\omega \rangle|} \geq {c \over {|k|^{\sigma_{0}}}}},\ \ \ \ c,\sigma_{0}> 0
\end{array}
\end{equation}
for all integer vectors $k\neq 0$.\ It is not difficult to show that for $\sigma_{0} >n$,\ the Lebesgue measure of the set of $\omega$ satisfying the above inequalities is positive for a suitably small $c$.

\subsection{The main result}

First we give the following definitions.

\begin{definition}\label{def2.4}
Let $\mathfrak{M}$ be a mapping given by (\ref{M}). It is said that $\mathfrak{M}$ has the intersection property if
$$\mathfrak{M}(\mathbf{\Gamma}) \cap \mathbf{\Gamma} \neq \emptyset$$
for every  curve $\mathbf{\Gamma}:\theta=\xi+\varphi(\xi),\ r=\psi(\xi)$,\ where the continuous functions $\varphi$ and $\psi$ are quasi-periodic in $\xi$ with the frequency $\omega=(\omega_1,\omega_2,\cdots,\omega_n)$.
\end{definition}

\begin{definition}\label{def2.400}
Let $\mathfrak{M}$ be a mapping given by (\ref{M}). We say that $\mathfrak{M}:\mathbb{R} \times [a,b]\to \mathbb{R}^2$  is an exact symplectic if
$\mathfrak{M}$ is  symplectic with respect to the usual symplectic structure $dr\wedge d\theta$ and for every  curve $\mathbf{\Gamma}:\theta=\xi+\varphi(\xi),\ r=\psi(\xi)$,\ where the continuous functions $\varphi$ and $\psi$ are quasi-periodic in $\xi$ with the frequency $\omega=(\omega_1,\omega_2,\cdots,\omega_n)$, we have
$$
\lim_{T\to+\infty}\frac{1}{2T}\int_{-T}^T\, rd\theta=\lim_{T\to+\infty}\frac{1}{2T}\int_{-T}^T\, r_1d\theta_1.
$$
\end{definition}

We claim that if the mapping $\mathfrak{M}$ is an exact symplectic map, then it has intersection property. In order to  prove this result,  we first give an useful  lemma, and its proof is simple.

\begin{lemma}\label{lem 7.1}
If $r = r(\theta)$ is quasi-periodic in $\theta$  and $F(\theta,r)$ is quasi-periodic in $\theta$ with the same frequency, then $F(\theta,r(\theta))$ is also quasi-periodic in $\theta$ with the same frequency.
\end{lemma}

Now we are going to prove the following lemma.

\begin{lemma}\label{lem 7.10}
If the mapping $\mathfrak{M}$ is an exact symplectic map, then it has intersection property.
\end{lemma}

\Proof Since the mapping $\mathfrak{M}$ is exact symplectic and it is also a monotonic twist map, according to the paper by Zharnitsky \cite{Zharnitsky00}, there is a function $H$ such that the mapping $\mathfrak{M}$ can be written by
$$r = -\frac{\partial }{\partial \theta} H(\theta_1-\theta, \theta),\quad r_1 = \frac{\partial}{\partial \theta_1}H(\theta_1-\theta,\theta),$$
where $H$ is quasi-periodic in the second variable.

Now we prove the intersection property of the mapping $\mathfrak{M}$, that is, given any continuous quasi-periodic curve $\Gamma : r=r(\theta)$, we need to prove that $\mathfrak{M}(\Gamma)\cap\Gamma\ne\emptyset.$  Define  two sets $\mathbb{B}$ and ${\mathbb{B}}_1$ : the set $\mathbb{B}$ is bounded by four curves $\big\{(\theta,r) : \theta=t\big\}$, $\big\{(\theta,r) : \theta=T\big\}$, $\big\{(\theta,r) : r=r_{*}\big\}$ and $\big\{(\theta,r) : r=r(\theta)\big\}$, the set ${\mathbb{B}}_1$ is bounded by four curves   $\big\{(\theta,r) : \theta=t\big\}$, $\big\{(\theta,r) : \theta=T\big\}$, $\big\{(\theta,r) : r=r_{*}\big\}$ and the image of $\Gamma$ under $\mathfrak{M}$. Here we choose $r_{*}<\min r(\theta)$. It is easy to show that the difference of the areas of ${\mathbb{B}}_1$ and $\mathbb{B}$ is
$$\Delta(t,T)=\int _{t}^{T}r_{1}d \theta_1-\int _{t}^{T}r d \theta=H(\theta_{1}(T)-T,T)-H(\theta_{1}(t)-t,t).$$
From the definition of $\mathfrak{M}$ and Lemma \ref{lem 7.1}, we know that $\theta_{1}(T)-T=r(T)+f(T,r(T))$ is quasi-periodic in $T$ and $\theta_{1}(t)-t=r(t)+f(t,r(t))$ is quasi-periodic in $t$. Hence using Lemma \ref{lem 7.1} again, it follows that  $\Delta(t,T)$ is quasi-periodic in $t$ and $T$.

Hence there are at least two pairs of $(t_1,T_1)$ and $(t_2,T_2)$ such that $\Delta(t_1,T_1)<0,\Delta(t_2,T_2)>0.$ The intersection property of $\mathfrak{M}$ follows from this fact, which proves the lemma. \qed

For the quasi-periodic mapping $\mathfrak{M}$  we assume that $f,g: \mathbb{R}^2 \mapsto \mathbb{R}$ are of class $\mathcal{C}^{p}$, and define
\begin{eqnarray*}
&&|x|=\max{(|\theta|,|r|)} \ \ \ \  \mbox{for}\ \ \   x=(\theta,r) \in \mathbb{R}^{2}, \\[0.2cm]
&&\big|h\big|_{\mathbb{R}^{2}}=\sup \limits _{x\in {\mathbb{R}^{2}}}|h(x)|, \\[0.2cm]
&&\|h\|_{p}=\sum\limits_{|k|\leq p}\,\sup \limits _{x\in \mathbb{R}^{2}}|D^{k}h(x)|
\end{eqnarray*}
if $p\geq 0$ is an integer, and
\begin{eqnarray*}
&&\|h\|_{p}=\sup \limits _{\substack{x\neq y\\ |k|=l}}{|D^{k}h(x)-D^{k}h(y)|\over |x-y|^{s}}+\sum\limits_{|k|\leq l}\,\sup \limits _{x\in \mathbb{R}^{2}}|D^{k}h(x)| \\[0.2cm]
\end{eqnarray*}
if $p=l+s$,\  $l \geq 0$ is an integer,\ $s \in (0,1) $,\ where
\begin{eqnarray*}
&&D^{k}={{\Big({\partial \over {\partial \theta}}\Big)}^{k_{1}}}\circ{{\Big({\partial \over {\partial r}}\Big)}^{k_{2}}},\ \ \ |k|=|k_{1}|+|k_{2}|,\ \ \ k=(k_{1},k_{2}).
\end{eqnarray*}

We choose a rotation number $\alpha$ satisfying the inequalities
\begin{equation}\label{b3}
\begin{array}{ll}
\left\{\begin{array}{ll}
a+12^{-3}\gamma\leq \alpha\leq b-12^{-3}\gamma,\\[0.4cm]
\Big|{\langle k,\omega \rangle {\alpha \over {2\pi}}-j}\Big|\geq {\gamma \over {|k|^\tau}},\ \ \ \ \mbox{for all}\ \  k \in \mathbb{Z}^n\backslash\{0\},\ \ \ \  j \in \mathbb{Z}
 \end{array}\right.
\end{array}
\end{equation}
with some constants $\gamma,\tau$ satisfying
\begin{equation}\label{b4}
\begin{array}{ll}
0<\gamma<{1\over 2} \min \{1,\ 12^{3}(b-a) \},\ \ \ \ \tau > n.
\end{array}
\end{equation}

Now we are in a position to state our main result.

\begin{theorem}\label{thm2.6}
Suppose that the quasi-periodic mapping $\mathfrak{M}$ given by (\ref{M}) is of class $\mathcal{C}^{p}\ (p>2\tau+1)$, and satisfies the intersection property, the functions $f(\theta,r),g(\theta,r)$ are quasi-periodic in $\theta$ with the frequency $\omega=(\omega_1,\omega_2,\cdots, \omega_n)$, and satisfy  the following smallness conditions
\begin{equation}\label{b5}
\begin{array}{ll}
|f|_{\mathbb{R}^{2}}+|g|_{\mathbb{R}^{2}}\leq  {6^{-(n+1)}\over 3}\ {q\over {300c_{0}}}{{\Big({1\over 72}\Big)}^{\tau}}{{\Big({\gamma\over \Gamma(\tau+1)}\Big)}^{2}},
\end{array}
\end{equation}
\begin{equation}\label{b6}
\begin{array}{ll}
\big\|f\big\|_{p}+\big\|g\big\|_{p}\leq  {6^{-(n+1)}\over 3}\ {{q(1-q)}\over {3600(3c_{1}+c_{2})}}{{\Big({1\over 288}\Big)^{\tau}}}{{\Big({\gamma\over \Gamma(\tau+1)}\Big)}^{2}},
\end{array}
\end{equation}
where $\Gamma$ is the  Gamma function, $\gamma,\tau$ satisfy (\ref{b4}), $c_{0},c_{1},c_{2}$ are positive constants depending only on $p$ and $\omega$, and $q$ is  a  number  satisfying
\begin{equation}\label{b7}
\begin{array}{ll}
0<q\leq \min \Big\{{{p-2\tau-1}\over {p+1}}\log 2,\ 10^{-2}4^{-\tau}\Big\}.
\end{array}
\end{equation}
Then for any number $\alpha$ satisfying the inequalities (\ref{b3}), the quasi-periodic mapping $\mathfrak{M}$ has an invariant curve $\Gamma_0$ with the form
$$
\begin{array}{ll}
\left\{\begin{array}{ll}
\theta=\theta'+\varphi(\theta'),\\[0.2cm]
r=\psi(\theta'),
 \end{array}\right.\ \ \
\end{array}
$$
where $\varphi, \psi$ are quasi-periodic  with the frequency $\omega=(\omega_1,\omega_2,\cdots,\omega_n)$, and the invariant curve $\Gamma_0$ is continuous and  quasi-periodic  with the frequency $\omega=(\omega_1, \omega_2, \cdots,\omega_n)$.\ Moreover,\ the  restriction of $\mathfrak{M}$ onto $\Gamma_0$ is
 $$\mathfrak{M}|_{{\Gamma} _{0}}: \theta_{1}^\prime=\theta^\prime +\alpha.$$
\end{theorem}

\begin{remark}
Here we assume that the mapping $\mathfrak{M}$ is of class $\mathcal{C}^{p}$ with $p>2\tau+1>2n+1$. $n=1$ corresponds to the periodic case, in which $p>3$ is the optimal smoothness condition. Hence our smoothness assumption for quasi-periodic mappings is optimal in this sense.
\end{remark}

\begin{remark}
If all conditions of Theorem \ref{thm2.6} hold,\ then the mapping $\mathfrak{M}$ has many invariant curves ${\Gamma _{0}}$, which can be labeled by the form $$\mathfrak{M}|_{\Gamma _{0}}: \theta_{1}^\prime = \theta' +\alpha$$
 of the restriction of $\mathfrak{M}$ onto ${\Gamma _{0}}.$\ In fact, given any $\alpha$ satisfying the inequalities (\ref{b3}), there exists an invariant curve ${\Gamma _{0}}$ of $\mathfrak{M}$ which is quasi-periodic  with the frequency  $\omega=(\omega_1,\omega_2,\cdots,\omega_n)$,\ and the restriction of $\mathfrak{M}$ onto ${\Gamma _{0}}$ has the form
$$\mathfrak{M}|_{\Gamma _{0}}: \theta_{1}^\prime=\theta' +\alpha.$$
The existence of such $\alpha$ can be found in Lemma \ref{Measure estimate}.
\end{remark}

The constants $c_{0},c_{1},c_{2}$ in the main result depend on  how well functions  of  class $\mathcal{C}^{p}$ can be approximated by analytic ones.
\begin{lemma}\label{lem2.9}
Let $h(\cdot,y) \in \mathcal{C}^{p} $ be a quasi-periodic function with the frequency  $\omega=(\omega_1,\omega_2,\cdots,\omega_n)$,\ then for  any $\delta>0$, there  exists a holomorphic  function $h_{\delta} :{\mathbb{C}}^2\mapsto \mathbb{C},$ $h_{\delta}(\cdot,y)\in Q_{r}(\omega)$, $h_{\delta}({\mathbb{R}}^2)\subseteq \mathbb{R}$ such  that  the  following  inequalities

\begin{equation*}
\begin{array}{ll}
\left\{\begin{array}{ll}
\big|h_{\delta}\big|_{E_\delta}\leq c_{0}\big|h\big|_{{\mathbb{R}}^2},\\[0.4cm]
{\big|h-h_{\delta}\big|_{{\mathbb{R}}^2}}\leq c_{1}{\big\|h\big\|_{p}}{{\delta}^p},\\[0.4cm]
{\big|h_{\delta}-h_{\delta^{'}}\big|_{E_\delta}}\leq c_{2}{\big\|h\big\|_{p}}{{\delta^{'}}^p}
 \end{array}\right.
\end{array}
\end{equation*}
hold for $0<\delta<{\delta^{'}}$, where
$$E_\delta=\{(x,y)\in {\mathbb{C}}^2: |\Im\ x|<\delta,\ \ |\Im\ y|<\delta\},\ |\cdot|_{E_\delta}=\sup\limits_{z\in E}|\cdot(z)|,$$
$c_{0},c_{1},c_{2}$ are positive  constants  only  depending on $p,\omega$.
\end{lemma}

The detail proof of Lemma \ref{lem2.9} is given in the Appendix. The proof of Lemma \ref{lem2.9} is similar to the periodic case.\ When  $h\in \mathcal{C}^{p}$ is a periodic function,\ there are some detail proofs  of  Lemma  \ref{lem2.9}  available  in  the  literature,  for example, see  Moser \cite[p. 528-529]{Moser66},\  R\"{u}ssmann \cite[p. 74-78]{Russmann70},\  Zehnder \cite[p. 110-113]{Zehnder75}.

\subsection{The measure estimate}

\begin{lemma} \label{Measure estimate}
If  $\tau > n$,\ then for suitable small $\gamma$,\ the set of $\alpha$ satisfying (\ref{b3}) has positive measure.
\end{lemma}

\noindent\textbf{Proof}: Choose some  $n$-dimensional frequency vector $\omega=(\omega_1,\omega_2,\cdots,\omega_n)$ satisfying (\ref{b2}) and let $\mathcal{D}_{\gamma,\tau}^{\omega}$ denote the set of all $\alpha\in\mathbb{R}$ satisfying (\ref{b3}) with the fixed $\gamma$ and $\tau$. Then $\mathcal{D}_{\gamma,\tau}^{\omega}$ is the complement of the open dense set $\mathcal{R}_{\gamma,\tau}^{\omega}$, where
\begin{eqnarray*}
\mathcal{R}_{\gamma,\tau}^{\omega}&=&\bigcup \limits _{\substack{ k \in \mathbb{Z}^n\backslash\{0\}\\  j \in \mathbb{Z}}}\mathcal{R}_{\omega,\gamma,\tau}^{k,j}\\&=&\bigcup \limits _{\substack{ k \in \mathbb{Z}^n\backslash\{0\}\\  j \in \mathbb{Z}}}\Big\{\alpha\in [a+12^{-3}\gamma,b-12^{-3}\gamma]:\big|{\langle k,\omega \rangle {\alpha \over {2\pi}}-j}\big|<{\gamma \over {|k|^\tau}}\Big\}.
\end{eqnarray*}

Now we estimate the measure of the set $\mathcal{R}_{\omega,\gamma,\tau}^{k,j}$. Set $|k_{\max}|=\max\limits_{1\leq i\leq n}|k_i|$, then there exists some  $1\leq m\leq n$ such that $|k_m|=|k_{\max}|$, and
$1\leq{{|k|}\over {|k_{\max}|}}\leq n.$  Therefore, we have
\begin{eqnarray*}
\mathcal{R}_{\omega,\gamma,\tau}^{k,j}&=&\Big\{\alpha\in [a+12^{-3}\gamma,b-12^{-3}\gamma]:\big|{\langle k,\omega \rangle {\alpha \over {2\pi}}-j}\big|<{\gamma \over {|k|^\tau}}\Big\}\\[0.2cm]
&=& \Big\{\alpha\in [a+12^{-3}\gamma,b-12^{-3}\gamma]:\big|k_{\max}\omega_m{\alpha \over {2\pi}}+\sum_{i\not =m} k_i\omega_i{\alpha \over {2\pi}}-j\big|<{\gamma \over {|k|^\tau}}\Big\}\\[0.2cm]
&=& \Big\{\alpha\in [a+12^{-3}\gamma,b-12^{-3}\gamma]:|k_{\max}||\omega_m||\alpha+b_{j}|<{2\pi\gamma \over {|k|^\tau}}\Big\}\\[0.2cm]
&=& \Big\{\alpha\in [a+12^{-3}\gamma,b-12^{-3}\gamma]:-b_{j}-\delta_k<\alpha<-b_{j}+\delta_k\Big\},
\end{eqnarray*}
where $b_{j}={1\over {k_{\max}|\omega_m|}}\Big\{\sum\limits_{i\not =m} k_i\omega_i\alpha-2\pi j\Big\}$ and
$
\delta_k = {2\pi\gamma \over {|k|^\tau}}\,{1\over {|k_{\max}||\omega_m|}}.
$
Hence,
$$\mbox{meas}\big(\mathcal{R}_{\omega,\gamma,\tau}^{k,j}\big)\le 2\delta_k = {4\pi\gamma \over {|k|^\tau}}\,{1\over {|k_{\max}||\omega_m|}}={4\pi\gamma \over {|k|^{\tau+1}}}\,{|k|\over {|k_{\max}|}}\,{1\over{|\omega_m|}}.
$$
Since $1\leq{{|k|}\over {k_{\max}}}\leq n,$ then we have the following measure estimate
$$\mbox{meas}\big(\mathcal{R}_{\omega,\gamma,\tau}^{k,j}\big)\leq O\bigg({\gamma \over {|k|^{\tau+1}}}\bigg).$$

Next we estimate the measure of the set $\mathcal{R}_{\gamma,\tau}^{\omega}$. Since for $\alpha\in\mathcal{R}_{\omega,\gamma,\tau}^{k,j}$,
$$\big|{\langle k,\omega \rangle {\alpha \over {2\pi}}-j}\big|<{\gamma \over {|k|^\tau}},$$
then we have
$$|j|\leq \big|\langle k,\omega \rangle \big|{{\alpha}\over {2\pi}}+{\gamma \over {|k|^\tau}}\leq c_0|k|,$$
where $c_0$ is a constant independent of $k$. Thus
\begin{eqnarray*}
\mbox{meas}(\mathcal{R}_{\gamma,\tau}^{\omega}) &\leq & \sum\limits_{k \in \mathbb{Z}^n\backslash\{0\}}\sum\limits_{\substack{j\in \mathbb{Z}\\ |j|\leq c_0|k|}}\mbox{meas}\big(\mathcal{R}_{\omega,\gamma,\tau}^{k,j}\big)\\[0.2cm]
&\leq& \sum\limits_{k \in \mathbb{Z}^n\backslash\{0\}}\sum\limits_{\substack{j\in \mathbb{Z}\\ |j|\leq c_0|k|}}O\bigg({\gamma \over {|k|^{\tau+1}}}\bigg)\leq \sum\limits_{k \in \mathbb{Z}^n\backslash\{0\}}O\bigg({\gamma \over {|k|^{\tau}}}\bigg).
\end{eqnarray*}
Also, if $\tau>n$,
$$\sum\limits_{k \in \mathbb{Z}^n\backslash\{0\}}{1\over {|k|^{\tau}}}\leq 2^n\sum_ {m=1}^{+\infty}{1\over {m^{\tau}}}\Big(\begin{matrix} n+m-1 \\ m \end{matrix} \Big)\leq 2^{2n-1}\sum_ {m=1}^{+\infty}{1\over {m^{\tau-n+1}}}<+\infty.$$
Hence, for any $\tau>n$,
$$\mbox{meas}(\mathcal{R}_{\gamma,\tau}^{\omega})\leq O(\gamma)$$
and
$$\mbox{meas}(\mathcal{D}_{\gamma,\tau}^{\omega})\to b-a \ \ \ \ \ \mbox{as} \ \ \ \ \ \gamma\to 0.$$
This completes the proof.\qed

\section{The iteration process}
In  this  section we present an iteration  process leading  to  the  proof of Theorem \ref{thm2.6}.

Firstly, we introduce  new variables by the  linear  transformation
$$
\begin{array}{ll}
\left\{\begin{array}{ll}
\theta=x,\\[0.1cm]
r=\alpha+{\varepsilon_{0}}y,
\end{array}\right.
\end{array}
$$
where $\alpha$ is the  chosen rotation  number satisfying  (\ref{b3}),\ ${\varepsilon_{0}}$ is  defined by
\begin{equation}\label{c2}
\begin{array}{ll}
{\varepsilon_{0}}=6^{-(\tau+{{n+1}\over 2})}\ {{\gamma}\over {{\Gamma(\tau+1)}}}.
\end{array}
\end{equation}

In  the  new coordinates  the  given mapping (\ref{M}) having the  intersection  property in  the  strip $S=\{(\theta,r)\in {\mathbb{R}}^2 : a<r<b \}$ gets the  form
$$
\begin{array}{ll}
A:\quad\left\{\begin{array}{ll}
x_1=x+\alpha+{\varepsilon_{0}}y+f(x,\alpha+{\varepsilon_{0}}y),\\[0.2cm]
y_1=y+{{\varepsilon_{0}}^{-1}}g(x,\alpha+{\varepsilon_{0}}y).
 \end{array}\right.
\end{array}
$$
Clearly the  intersection  property  is  preserved and holds  in  the  strip
\begin{equation}\label{c4}
\begin{array}{ll}
{S^{\star}}=\{(x,y)\in {\mathbb{R}}^2:|y|<600^{-1}\},
\end{array}
\end{equation}
where we have used (\ref{c2}) and ${\Gamma(\tau+1)}\geq  1$ for $\tau>n$.

Since the $\mathcal{C}^{p} $ functions $f(\cdot,\alpha+{\varepsilon_{0}}y),g(\cdot,\alpha+{\varepsilon_{0}}y)$ are quasi-periodic with the frequency  $\omega=(\omega_1,\omega_2,\cdots,\omega_n)$,\ by assumption we may apply Lemma  \ref{lem2.9} to obtain  a family  of holomorphic functions $f_{\delta}(\cdot,\alpha+{\varepsilon_{0}}y),\ g_{\delta}(\cdot,\alpha+{\varepsilon_{0}}y)\in Q_{r}(\omega)\ (\delta>0) $, with which we define the quasi-periodic mappings
$$
\begin{array}{ll}
{A_{\delta}}:\quad\left\{\begin{array}{ll}
x_1=x+\alpha+{\varepsilon_{0}}y+{f_{\delta}}(x,\alpha+{\varepsilon_{0}}y),\\[0.2cm]
y_1=y+{{\varepsilon_{0}}^{-1}}{g_{\delta}}(x,\alpha+{\varepsilon_{0}}y).
 \end{array}\right.
\end{array}
$$
Define
$$E_{\delta}=\{(x,y)\in {\mathbb{C}}^2  : |\Im\ x|<{\delta},\ |\Im\ y|<{\delta} \}$$
and a sequence
\begin{equation}\label{c6}
\begin{array}{ll}
{\delta_{k}}={\Big({{1+q}\over 2}\Big)}^{k},\ \ \ \ k=0,1,2,\cdots,
\end{array}
\end{equation}
where $q$ is  a real  number  satisfying (\ref{b7}),\ and set
$$
{E_{k}}={E_{{\delta_{k}}}},\ \ \ \  {A_{k}}={A_{{\delta_{k}}}},\ \ \ \   k=0,1,\cdots .
$$
Then the  estimates  of Lemma  \ref{lem2.9} can be written in  the  form
\begin{equation}\label{c7}
\begin{array}{ll}
\left\{\begin{array}{ll}
\big|A_{0}-\Omega_{0}\big|_{E_{0}}\leq {{\varepsilon_{0}}^{-1}}{c_{0}}\big(\big|f\big|_{{\mathbb{R}}^2}+\big|g\big|_{{\mathbb{R}}^2}\big),\\[0.4cm]
\big|A-A_{k}\big|_{{{\mathbb{R}}^2}}\leq {{\varepsilon_{0}}^{-1}}{c_{1}}\big(\big\|f\big\|_{p}+\big\|g\big\|_{p}\big){{\delta_{k}}^p},\\[0.4cm]
\big|A_{k}-A_{k+1}\big|_{E_{k+1}}\leq {{\varepsilon_{0}}^{-1}}{c_{2}}\big(\big\|f\big\|_{p}+\big\|g\big\|_{p}\big){{\delta_{k}}^p},
 \end{array}\right.\   k=0,1,2,\cdots ,
\end{array}
\end{equation}
where the  mapping
$$
\begin{array}{ll}
\Omega_{0}: \left\{\begin{array}{ll}
x_1=x+\alpha+{\varepsilon_{0}}y,\\[0.2cm]
y_1=y.
 \end{array}\right.
\end{array}
$$

Before we describe the  iteration  process,  some  definitions  and  notations  are useful.

\medskip

(i) Given subsets $D_{1},\cdots,D_{\ell}$ of ${\mathbb{C}^m}$ and functions $F_{j} : D_{j} \mapsto {\mathbb{C}^m},\ j=1,\cdots,\ell-1$.\ Then
\[ D_{1}\xrightarrow {F_{1}}     D_{2} \xrightarrow {F_{2}}     D_{3} \xrightarrow {F_{3}}     \cdots \xrightarrow {F_{\ell-1}}     D_{\ell}{} \]
exists  if
$$F_{j}(D_{j})\subseteq D_{j+1},\ j=1,\cdots, \ell-1.$$
In  the  case $F_{j}=\mbox{id}=\mbox{id}_{D_{j}}$ this  condition  means
$$D_{j}\subseteq D_{j+1},\ \ \ \mbox{id}(x)=x\ \ \ \mbox{for  all}\ \ \  x\in D_{j}.$$

\medskip

(ii) For $r,\ s> 0$,  define
$$D(r,s)=\{(x,y)\in {\mathbb{C}}^2: |\Im\ x|<r,\ |y|<s \}.$$

\medskip

(iii) For $D=D(r,s)$, denote by $T(D)=T(r,s)$ the set of all holomorphic functions $F : D \mapsto {\mathbb{C}^2}$ satisfying  the  identity
$$\sigma\circ F=F\circ \sigma\big|_{D},$$
where $\sigma$ is defined  by
$$(x,y)\mapsto \sigma(x,y)=(\bar{x},\bar{y})$$
for all $(x,y)\in {\mathbb{C}^2}$ with $x=a+bi$,  $\bar{x}=a-bi, a,b\in \mathbb{R}$.

\medskip

(iv) Define the  mappings $\Omega_{k}\ (k=0,1,\cdots)$ by
$$
\begin{array}{ll}
\Omega_{k}:\left\{\begin{array}{ll}
x_1=x+\alpha+{\varepsilon_{k}}y,\\[0.2cm]
y_1=y,
 \end{array}\right.\ \ \ \ \ (x,y)\in {\mathbb{C}^2},\ \ \ {\varepsilon_{k}}={2^{-k\tau}}{\varepsilon_{0}}.
\end{array}
$$

\medskip

(v) In ${\mathbb{C}^m}\ (m=1,2,\cdots)$, define the  norm
$$|x|=\max \limits_{j}|x_{j}|\ \ \ \mbox{ for} \ \ \ x=(x_{1},\cdots,x_{m})\in {\mathbb{C}^m}.$$

\medskip

(vi) Given $f : D \mapsto {\mathbb{C}^m} $ with $D=D(r,s)$,\  define
\begin{eqnarray*}
\big|f\big|_{\rho,\tilde{\sigma}}&=&\sup \limits _{\substack{|\Im\ x|< \rho\\ |y|<\tilde{\sigma}}}{\big|f(x,y)\big|},\ \ \ 0<\rho\leq r,\ \ \ 0<\tilde{\sigma}\leq s,\\
\big|f\big|_{D}&=&\big|f\big|_{r,s}.
\end{eqnarray*}

Now we are going back to  the  quasi-periodic mappings $A_{k} : E_{k}\mapsto {\mathbb{C}^2}$ defined above. We try  to  fix  domains
$$D_{k}=D(r_{k},s_{k}),\ \ \ D'_{k}=D(r'_{k},s'_{k}),$$
and to  find  mappings  $Z_{k}\in T(D'_{k}),H_{k}\in T(D_{k})$, and $Z_{k}-\mbox{id},H_{k}-\Omega_k$ are quasi-periodic with the frequency $\omega=(\omega_1,\omega_2,\cdots,\omega_n)$ in the first variable,   such that  the  diagrams
$$\begin{CD}
D_{k} @>{Z_{k}|_ {D_{k}}}>> E_{k}\\
@VH_{k}VV  @VVA_{k}V\\
{D'_{k}}@>>Z_{k}> {\mathbb{C}^2}
\end{CD} \eqno (3.5)_{k}$$
exist and commute for $k=0,1,\cdots .$

A proper choice  for  the  constants $r_{k},\ s_{k},\ r'_{k},\ s'_{k}$ is
$$
\left\{\begin{array}{ll}
r_{k}=2^{-k},\ \ s_{k}=2^{-k}s_{0},\ \ s_{0}=300^{-1}2^{-\tau},\\[0.2cm]
r'_{k}={4\over 3}(r_{k}-s_{k}),\ \ s'_{k}={4\over 3}s_{k},\ \ k=0,1,\cdots .
 \end{array}\right. \eqno (3.6)$$
Then obviously $D_{k}\subseteq D'_{k}\ (k=0,1,\cdots).$  About the  mappings $Z_{k},\ H_{k}$ the  following relations  are needed
$$Z_{k}(D_{k})\subseteq D_{0},  \eqno (3.7)_{k}$$
$$
b_{k}|\zeta-\zeta'|\leq |Z_{k}(\zeta)-Z_{k}(\zeta')|\leq B_{k}|\zeta-\zeta'|,\quad  \ \ \zeta,\ \zeta' \in D'_{k},
\eqno (3.8)_{k}
$$

$$
|Z_{k+1}(\zeta)-Z_{k}(\zeta)|\leq {2\over 3}q B_{k}s_{k},\quad  \ \  \zeta=(\xi,0) \in D_{k+1},
\eqno (3.9)_{k}
$$
$$
\Big|H_{k}-\Omega_{k}\Big|_{D_{k}}\leq M_{k},
\eqno (3.10)_{k}
$$
where
$b_{k}=2^{-k\tau}(1-q)^{k}$, $B_{k}=(1+q)^k$, $M_{k}=2^{-k(\tau+1)}M_{0},\ M_{0}={1\over 3}q \varepsilon_{0} s_{0}.$

Finally define $Z_{0}=\mbox{id}_{D'_{0}}\in T(D'_{0})$, then by means of this  iteration  process,\ if  it  exists,\ the  assertion of Theorem
\ref{thm2.6} can easily be proved.

In  fact,  from (3.6), $(3.8)_{k}$ and $(3.9)_{k}$, the  sequence $Z_{0},Z_{1},\cdots$ converges uniformly on $\mathbb{R}\times \{0\}$ and the  limit
$Z_{\infty}(\xi)=\lim \limits_{k\rightarrow\infty}Z_{k}(\xi,0)$
is  continuous on $\mathbb{R}$.

Since $Z_{k}\in T(D'_{k})\ (k=0,1,\cdots) $,
$$\sigma\circ Z_{\infty}=Z_{\infty}\circ \sigma \big|_ {\mathbb{R}},\ \ \ \ Z_{\infty}(\mathbb{R})\subseteq \mathbb{R}^2.$$
Now the  commutativity  of $(3.5)_{k}$ yields
$$Z_{k}\circ H_{k}=A_{k}\circ Z_{k}\big|_ {D_{k}}.$$
Hence
$$A\circ Z_{k}(\xi,0)-Z_{k}(\Omega_{k}(\xi,0))=(A-A_{k})(Z_{k}(\xi,0))+Z_{k} \big(H_{k}(\xi,0)\big)-Z_{k} \big(\Omega_{k}(\xi,0)\big)$$
and by virtue  of $(3.8)_{k},(3.10)_{k}$ consequently
$$\big|A\circ Z_{k}(\xi,0)-Z_{k}(\alpha+\xi,0)\big|\leq \big|A-A_{k}\big|_{\mathbb{R}^2}+B_{k}M_{k},\ \ \ \xi\in \mathbb{R}.$$
Passing to  the  limit  we get
$$A\circ Z_{\infty}(\xi)=Z_{\infty}(\alpha+\xi),\ \ \ \xi\in \mathbb{R}$$
in  view of (\ref{c6}),\ (\ref{c7}) and $B_{k}M_{k}\rightarrow 0.$

Therefore,\ we can obtain the existence of invariant curves of the mapping $A$,\ and from the relationship between the mappings $A$ and $\mathfrak{M}$,  one  can also get  the existence of invariant curves of the mapping $\mathfrak{M}$.

From the above analysis, firstly we need to  prove the  assertion
$$
\left\{\begin{array}{ll}
\text{The\ diagram}\ (3.5)_{k}\  \text{exists\  and\  commutes}\\[0.4cm]
\text{with\ some}\ Z_{k}\in T(D'_{k}),\ H_{k}\in T(D_{k})\\[0.4cm]
\text{satisfying}\ (3.7)_{k},\ (3.8)_{k},\ (3.10)_{k}
 \end{array}\right. \eqno (3.11)_{k}$$
and the estimate $(3.9)_{k}$ for $k=0,1,\cdots$.

The proofs of $(3.9)_{k}$ and $(3.11)_{k}\ (k=0,1,\cdots)$ are  done by the complete  induction. Let us first  consider the  case $(3.11)_{0}.$ As a consequence  of the  definition of $Z_{0}$, the
relations $(3.7)_{0}$ and $(3.8)_{0}$ are obvious.  Moreover if  we define
$H_{0}=A_{0}\big|_{ D_{0}},$\ by virtue  of $(\ref{b5}),\ (\ref{c2}),\ (\ref{c7}),\ (3.6)$, then
$$\big|A_{0}-\Omega_{0}\big|_{E_{0}}\leq {{\varepsilon_{0}}^{-1}}{c_{0}}\big(\big|f\big|_{{\mathbb{R}}^2}+\big|g\big|_{{\mathbb{R}}^2}\big)=M_{0},$$
which is the wanted  estimate $(3.10)_{0}$.\

From the definition of $\varepsilon_{0}$ and $q$,\ we have
$$M_{0}={1\over 3}q\varepsilon_{0}s_{0}<{1\over 3}s_{0}.$$
Hence
$$H_{0}(D_{0}) \subseteq  D(r_{0}+\varepsilon_{0}s_{0}+M_{0},s_{0}+M_{0})
               \subseteq D(1+s_{0}+{1\over 3}s_{0},s_{0}+{1\over 3}s_{0})\subseteq D'_{0}.   $$
From the definitions of $D_{0}$ and $E_{0}$, $D_{0}\subseteq E_{0}$ holds.\ Therefore
$$
D_{0}\subseteq E_{0},\ \ \ \ H_{0}(D_{0})\subseteq {D'_{0}}.
$$
Thus,\ the  diagram  $(3.5)_{0}$ exists and commutes.

Now let  us suppose that  $(3.11)_{k}$ is  true  for  some  $k\geq 0.$ We have to  show $(3.11)_{k+1}$ and $(3.9)_{k}.$ On this  way the  crucial  result is  the  construction  of the  commuting diagram
 $$\ \ \ \ \ \ \ \  \ \ \  \ \ \ \  \ \ \ \ \ \ \ \ \ \ \ \ \ \ \ \ \ \ \ \ \ \ \ \ \ \ \ \ \ \ \ \ \xymatrix{
 D_{k+1} \ar [rr]^{W_{k}\big|_ { D_{k+1}}} \ar "2,2"_{\Phi_{k+1}}&{}&
 D_{k}\ar [dd]^{H_{k}}\\
 {}&{D^\star_{k+1}}\ar "3,1"_{\mbox{id}}&{}&\ \ \ \ \ \ \ \ \ \ \ \ \ \ \ \ \ \ \ \     {(3.12)}\\
{D'_{k+1} } \ar [rr]_{W_{k}}&{}
&{D'_{k}}}
$$
with
$$
D^\star_{k+1}=D(r'_{k+1}-{1\over 7}s_{k},s'_{k+1}-{1\over 7}s_{k})
$$
and mappings $W_{k}\in T(D'_{k+1}),\  \Phi_{k+1}\in T(D_{k+1})$, and $W_{k}-\mbox{id},\Phi_{k+1}-\Omega_{k+1}$ are quasi-periodic with the frequency $\omega=(\omega_1,\omega_2,\cdots,\omega_n)$ in the first variable.

The existence of  the commuting diagram (3.12) is  guaranteed  by the inductive theorem (Theorem \ref{thm5.3}),
which we will prove in  Section 5.\ This theorem  also gives the  following  estimates
$$
{2^{-\tau}}(1-q)|\zeta-\zeta'|\leq |W_{k}(\zeta)-W_{k}(\zeta')|\leq (1+q)|\zeta-\zeta'|,\ \   \zeta,\ \zeta' \in D'_{k+1},
\eqno (3.13)$$
$$
 |W_{k}(\zeta)-\zeta|\leq {2\over 3}q s_{k},\quad    \zeta=(\xi,0) \in D_{k+1},
 \eqno (3.14)$$
$$
 \big|\Phi_{k+1}-\Omega_{k+1}-Q_{k}\big|_{D_{k+1}}\leq {5\over 24}M_{k+1},
\eqno (3.15)$$
where $Q_{k}$ is  a polynomial  of degree 2 in  the  second variable only
$$
 Q_{k}(\eta)=(0,a_{0k}+a_{1k}\eta +a_{2k}{\eta^2}),\ \ a_{0k},a_{1k},a_{2k}\in \mathbb{R}.
\eqno (3.16)$$

With these  assertions of the inductive theorem  we can  show $(3.11)_{k+1}$ and $(3.9)_{k}.$  From the  diagrams $(3.5)_{k}$ and (3.12) we see that

$$\xymatrix{
D_{k+1} \ar [rr]^{W_{k}\big|_ {D_{k+1}}} \ar "2,2"_{\Phi_{k+1}}&{}&
 D_{k}\ar [dd]^{H_{k}} \ar [rr]^{Z_{k}} &{}& E_{k} \ar [dd]^ {A_{k}}\\
{}&{D^\star_{k+1}}\ar "3,1"_{\mbox{id}}&{}&\\
 {D'_{k+1} }\ar [rr]_{W_{k}}&{}
&{D'_{k}}\ar [rr]_{Z_{k}}&{}
&{\mathbb{C}^2}.
}$$
Define
$$Z_{k+1}=Z_{k}\circ W_{k},$$
of course, $Z_{k+1}\in T(D'_{k+1})$, $Z_{k+1}-\mbox{id}$ is quasi-periodic with the frequency $\omega=(\omega_1,\omega_2,\cdots,\omega_n)$ in the first variable, and $(3.7)_{k+1}$ holds.  The proof of $(3.8)_{k+1}$ by means of $(3.8)_{k}$ and (3.13) is  obvious  if  we notice
$W_{k}(D'_{k+1})\subseteq D'_{k}$ in (3.12).\ The inequality  $(3.9)_{k}$ follows  from $(3.8)_{k}$,\ $(3.14)$ and $W_{k}(D_{k+1})\subseteq D'_{k}$ as a consequence of $D_{k+1}\subseteq D'_{k+1}.$ We also need to prove
$$Z_{k+1}(D_{k+1})\subseteq E_{k+1}.$$
In fact,\
$$\big|\Im\ Z_{k+1}(\zeta)\big|=\big|\Im\ Z_{k+1}(\Re\ \zeta+i\, \Im\ \zeta)\big| ,\ \ \ \zeta\in D_{k+1},$$
where $Z_{k+1}$ is real for real arguments as an element of $T'(D_{k+1})$,\ we get
$$\big|\Im\ Z_{k+1}(\Re\ \zeta+i\, \Im\ \zeta)\big|=\big|\Im\ \big\{Z_{k+1}(\Re\ \zeta+i \,\Im\ \zeta)-Z_{k+1}(\Re\ \zeta)\big\}\big|,$$
hence
$$\big|\Im\ Z_{k+1}(\zeta)\big|\leq \big| Z_{k+1}(\zeta)- Z_{k+1}(\Re\ \zeta)\big|.$$
Using $(\ref{c6})$,\ $(3.6)$,\ $(3.8)_{k}$ we obtain
\begin{eqnarray*}
\big|\Im\ Z_{k+1}(\zeta)\big| &\leq& \big| Z_{k+1}(\zeta)- Z_{k+1}(\Re\ \zeta)\big| \leq B_{k+1}\big|\Im\ \zeta\big|\\[0.2cm]
&\leq& B_{k+1}\max\big\{r_{k+1},s_{k+1}\big\}=\delta_{k+1}.
\end{eqnarray*}
By the definition of $E_{k}\ (k=0,1,2,\cdots)$, $Z_{k+1}(D_{k+1})\subseteq E_{k+1}.$

With these  assertions we obtain the  following commuting diagram
$$\ \ \ \ \ \ \ \ \ \ \ \ \ \ \ \ \ \ \ \ \ \ \ \ \ \ \ \ \ \ \ \ \ \ \ \ \xymatrix{
D_{k+1} \ar [rr]^{Z_{k+1}\big|_ {D_{k+1}}} \ar "2,2"_{\Phi_{k+1}}&{}&
 E_{k+1}\ar [dd]^{{A_{k}\big|_ {E_{k+1}}}}\\
{}&{D^\star_{k+1}}\ar "3,1"_{\mbox{id}}&{}&\ \ \ \ \ \ \ \ \ \  \ \ \ \ \ \ \ \ \ \ \ {(3.17)}\\
 {D'_{k+1} }\ar [rr]_{Z_{k+1}}&{}
&{\mathbb{C}^2}.
}$$

Comparing the  diagrams $(3.5)_{k+1}$ and (3.17),\ it  remains  to  prove that  we can replace $\Phi_{k+1}$ by $H_{k+1}$ if  we replace $A_{k}\big|_ {E_{k+1}}$ by $A_{k+1}.$ Moreover we have to  show $(3.10)_{k+1}$, which is not possible without going back to  the  original quasi-periodic mapping $A$ in  order to  use the intersection  property,  and to  estimate the  polynomial  (3.16)  well enough such that $(3.10)_{k+1}$ follows  from  (3.15).

In the  following,\ we will prove these  assertions. First of all  we give some useful definitions and lemmas.\ $\mathbb{F}$ stands for $\mathbb{C}$ or $\mathbb{R}$,\ we call a function analytic
if  it  is  holomorphic in  the case $\mathbb{F}=\mathbb{C}$ and $\mathbb{R}$-analytic in the case $\mathbb{F}=\mathbb{R}.$
Moreover, for $D\subseteq \mathbb{F}^m$ and $d>0$, define the set
$$D-d=\Big\{x\in \mathbb{F}^m:  \big\{ y\ :\ |y-x|\leq d\big\} \subseteq D \Big\},$$
which may be empty. If $D$  is  open so is $D-d.$

\begin{lemma} [Lemma 2 in {\cite{Russmann83}}]\label{lem3.1}
Let $D$ be an open subset of $\mathbb{F}^m$,\ and $F : D\mapsto F(D)\subseteq \mathbb{F}^m$ be an analytic mapping satisfying the  estimate
$$
b|z-z'|\leq |F(z)-F(z')|
$$
for all $z,\ z'\in D $ and some $b>0.$ Then $F(D)$ is  open, and the  inverse  mapping $F^{-1}:F(D) \mapsto D$ exists and is  analytic.  Moreover for  any $d>0$ we have
$$
F(D-d)\subseteq F(D)-bd.
$$
\end{lemma}

It  is  useful  to  introduce  an arbitrary bijection
$\Lambda : \mathbb{F}^m \rightarrow \mathbb{F}^m$
and to  denote by
$\Delta=\Delta(\Lambda,\mathbb{F}^m)$
the  set of all  subsets $D$ of $\mathbb{F}^m$ which are invariant  under $\Lambda.$  Furthermore denote by
$\Sigma=\Sigma(\Lambda,\mathbb{F}^m)$
the  class of all  functions $F:D\mapsto \mathbb{F}^m$ such that
$D\in \Delta,\ \Lambda\circ F=F\circ \Lambda \big|_ D.$
Clearly, $F(D)\in\Delta$ for  a function  of class $\Sigma$  which domain is $D,$ and  if  $F$ is  injective, then  also $F^{-1}: F(D)\mapsto\mathbb{F}^m$ is of class $\Sigma.$ Moreover if $F:D\mapsto \mathbb{F}^m,$\ $G:E\mapsto \mathbb{F}^m$ are of class $\Sigma,$ and if $G(E)\subseteq D$, then $F\circ G$ is of class $\Sigma.$

\begin{lemma} [Theorem 3 in {\cite{Russmann83}}]\label{lem3.2}
Let $D,\ D',\ E$ be open subsets of $\mathbb{F}^m$ belonging to $\Delta$ with $D\subseteq D'$ and $A',\ Z,\ \Phi'$ be analytic mappings of class $\Sigma$ such that  the  diagram
$$\ \ \ \ \ \ \ \ \ \ \ \ \ \ \ \ \ \ \ \ \ \ \ \ \ \ \ \ \ \ \ \ \ \ \  \ \ \ \ \  \xymatrix{
D \ar [rr]^{Z\big|_ D} \ar "2,2"_{\Phi'}&{}&
 E\ar [dd]^{A'}\\
{}&{D'-d}\ar "3,1"_{\mbox{id}}&{}&\ \ \ \ \ \ \ \ \ \ \ \  \ \ \ \  \ \ \ \ \  {(3.18)}\\
 {D' }\ar [rr]_{Z}&{}
&{\mathbb{F}^m}
}$$
exists and commutes  with some $d>0,$ and the  estimate
$$b|\zeta-\zeta'|\leq |Z(\zeta)-Z(\zeta')|$$
holds for  all $\zeta,\zeta'\in D'$ with some $b>0.$ Then for  any continuous mapping $A'':E \mapsto \mathbb{F}^m$ of class $\Sigma$ satisfying  the  estimate
$$\big|A'-A''\big|_{E}\leq bd,\eqno (3.19)$$
there  exists a continuous  mapping $\Phi''$   of a class $\Sigma$ such  that the diagram
$$\begin{CD}
D @>{Z\big|_ D}>> E\\
@V\Phi''VV  @VVA''V\\
D'@>>Z> \mathbb{F}^m
\end{CD}$$
exists and commutes,  and the  estimate
$$\big|\Phi'-\Phi''\big|_{D}\leq{b^{-1}}\big|A'-A''\big|_{E}\eqno (3.20)$$
is  valid. If  $A''$ is  analytic so  is $\Phi''.$
\end{lemma}

We apply Lemma \ref{lem3.2} to  the  diagram (3.17) in  order to  obtain $(3.5)_{k+1}.$ Set
$$
A'=A_{k}\big|_{ E_{k+1}},\ \ A''=A_{k+1},\ \ \Phi'=\Phi_{k+1},\ \ Z=Z_{k+1},\eqno (3.21)$$
$$D=D_{k+1},\ \ D'=D'_{k+1},\ \ E=E_{k+1},$$
$$m=2,\ \ \mathbb{F}=\mathbb{C},\ \ b=b_{k+1},\ \ d={1\over 7}s_{k}.
$$
Moreover, $\Delta$ is the  set of all  subsets of $\mathbb{C}^2$ which are invariant under $\Lambda=\sigma$ such that $\Sigma$ is  the class  of all  functions $F:D\mapsto \mathbb{C}^2$ with $D\in\Delta$ and $\sigma\circ F=F\circ \sigma.$ Then $D_{k+1},\  D'_{k+1},\ E_{k+1}$ are open sets belonging to $\Delta,$ and (3.21) represents  analytic functions  of class $\Sigma.$ For $A'$ and $A''$ this follows  from  Lemma  \ref{lem2.9},  for $\Phi'$ and $Z$ this is true  because of $\Phi_{k+1}\in T(D_{k+1}),Z_{k+1}\in T(D'_{k+1}).$ In  addition $(3.8)_{k+1}$ is  valid.\  Therefore  Lemma  \ref{lem3.2} can be applied and gives a function $H_{k+1}=\Phi''\in T(D_{k+1})$, $H_{k+1}-\Omega_{k+1}$ is quasi-periodic with the frequency $\omega=(\omega_1,\omega_2,\cdots,\omega_n)$ in the first variable,  such  that  the diagram $(3.5)_{k+1}$ exists and commutes.

We also need to prove $H_{k+1}$ satisfies $(3.10)_{k+1}$. By  $(\ref{b6}),\ (\ref{c2}),\ (\ref{c6}),\ (\ref{c7})$,\ we obtain
$$\big|A_{k+1}-A_{k}\big|_{E_{k+1}}\leq {1\over 3} \varepsilon_{0} c_{2} {{q(1-q)}\over {600(3c_{1}+c_{2})}}{{(1+q)^{kp}}}{2^{-3\tau-kp}}.$$
The  necessary condition  which we have to  require  is
$$
{{(1+q)^p}\over {1-q}}\leq 2^{p-1-2\tau},
\eqno (3.22)$$
for then  we get  the estimate
$$\big|A_{k+1}-A_{k}\big|_{E_{k+1}}\leq {1\over 3} \varepsilon_{0} c_{2} {{q{(1-q)^{k+1}}}\over {600(18c_{1}+6c_{2})}}{2^{-3\tau-k-2k\tau}}.$$
By  the  definitions  of $s_{0},\ b_{k+1},\ M_{k+1}$ in (3.6), $(3.8)_{k+1}$ and $(3.10)_{k+1},$ we have
$$b_{k+1}=2^{-(k+1)\tau}(1-q)^{(k+1)},\  M_{k+1}={1\over 3}q \varepsilon_{0}{300^{-1}}{2^{-\tau}}2^{-(k+1)(\tau+1)},$$
hence
$$\big|A_{k+1}-A_{k}\big|_{E_{k+1}}\leq{{ c_{2}}\over {18c_{1}+6c_{2}}}{b_{k+1}}M_{k+1}.
\eqno (3.23)$$
Since $M_{k+1}s_{k}^{-1}\leq M_{0}s_{0}^{-1}\leq{3^{-1}}$,\ and put $d={1\over 7}s_{k},$\ we get
$$\big|A_{k+1}-A_{k}\big|_{E_{k+1}}\leq b_{k+1}d.$$
According to Lemma \ref{lem3.2},\ we have
$$
\big|H_{k+1}-\Phi_{k+1}\big|_{D_{k+1}}\leq {{c_{2}}\over {18c_{1}+6c_{2}}}M_{k+1}.
$$
By $(3.15)$,\ we get
$$
 \big|\Phi_{k+1}-\Omega_{k+1}-Q_{k}\big|_{D_{k+1}}\leq {5\over 24}M_{k+1},$$
then we have the estimate
$$
 \big|H_{k+1}-\Omega_{k+1}-Q_{k}\big|_{D_{k+1}}\leq \Big({{5\over 24}+{{c_{2}}\over {18c_{1}+6c_{2}}}}\Big)M_{k+1}
$$
with a polynomial $Q_{k}$ defined in  (3.16).

In  order to  obtain a proper estimate for $Q_{k}$,\ we apply Lemma \ref{lem3.2} once more to  the diagram (3.17), where this  time  we restrict $D_{k+1}$ to $D=\mathbb{R}^2\cap D_{k+1}$ such  that  we consider (3.18) with
$$
\begin{array}{ll}
m=2,\ \ \mathbb{F}=\mathbb{R},\ \ b=b_{k+1},\ \ d={1\over 7}s_{k}, \\[0.2cm]
D=\mathbb{R}^2\cap D_{k+1},\ \ D'=\mathbb{R}^2\cap D'_{k+1},\ \ E=\mathbb{R}^2\cap E_{k+1}=\mathbb{R}^2,\\[0.2cm]
A'=A_{k}\big|_ {\mathbb{R}^2},\ \ \Phi'=\Phi_{k+1}\big|_{\mathbb{R}^2\cap D_{k+1}},\ \ Z=Z_{k+1}\big|_ {\mathbb{R}^2\cap D'_{k+1}}.
\end{array}
$$
Here we put $\Lambda=\sigma\big |_{\mathbb{R}^2}$ such that $D,\ D',\ E$ are open subsets of $\mathbb{R}^2$ belonging to $\Delta,$ and $A',\ \Phi,\ Z$ are analytic  functions  of class $\Sigma.$ Also the  original  quasi-periodic mapping  $A$ defined at the  beginning of Section 3 is  of class $\Sigma,$ and it  is  continuous. Therefore using $(3.7)_{k+1}$,\ Lemma \ref{lem3.2} is  again applicable  and we obtain a continuous  function $\Psi_{k+1}=\Phi''$ of class $\Sigma$ such that  the  diagram
$$\begin{CD}
\mathbb{R}^2\cap D_{k+1} @>{Z_{k+1}\big | _{\mathbb{R}^2\cap D_{k+1}}}>> \mathbb{R}^2\\
@V\Psi_{k+1}VV  @VVA_kV\\
\mathbb{R}^2\cap D'_{k+1}@>>{Z_{k+1}\big | _{\mathbb{R}^2\cap D'_{k+1}}}> \mathbb{R}^2
\end{CD}\eqno (3.24)$$
exists and commutes provided (3.19) can be satisfied. Furthermore by means of (\ref{c7}),\ (3.20) and (3.23) we get  the estimate
$$
\big|\Psi_{k+1}-\Phi_{k+1}\big|_{\mathbb{R}^2\cap D_{k+1}}\leq {{c_{1}}\over {18c_{1}+6c_{2}}}M_{k+1},$$
which leads  with (3.15) to
$$
 \big|\Psi_{k+1}-\Omega_{k+1}-Q_{k}\big|_{\mathbb{R}^2\cap D_{k+1}}\leq \Big({{5\over 24}+{{c_{1}}\over {18c_{1}+6c_{2}}}}\Big)M_{k+1}.\eqno (3.25)$$

We recall  that  the quasi-periodic mapping  $A$ has the  intersection  property at least  in  the  strip (\ref{c4}). We apply this  property to  the  family  of curves
$$\mathbb{R}\ni\ \ \xi \mapsto \ Z_{k+1}(\xi,\eta),\ \ \ \ \ -s_{k+1}<\eta<s_{k+1},$$
where it  is  clear that  these  curves lie  in $S^\star.$ Moreover these  curves
satisfy  the  conditions of Definition  \ref{def2.4}. For  each $\eta$ with $-s_{k+1}<\eta<s_{k+1}$,\ there  are  real  numbers
$\xi_{0},\ \xi_{1}$ such that
$$A\circ Z_{k+1}(\xi_{0},\eta)=Z_{k+1}(\xi_{1},\eta).$$
On the  other hand from  the  commuting diagram (3.24) we have
$$A\circ Z_{k+1}(\xi_{0},\eta)=Z_{k+1}(\Psi_{k+1}(\xi_{0},\eta)).$$
The  mapping $Z_{k+1}$ is  analytic, and  as  a
consequence of $(3.8)_{k+1}$  it  is  injective. Thus the  injectivity  of $Z_{k+1}$ yields
$$(\xi_{1},\eta)=\Psi_{k+1}(\xi_{0},\eta),$$
hence
$$\eta=\Psi^{(2)}_{k+1}(\xi_{0},\eta),\ \ \ \  -s_{k+1}<\eta<s_{k+1},$$
where (2)  indicates  the  second component of a vector.

This equation leads  to  a reasonable estimate for  the  polynomial (3.16).  Since we use the maximum norm we get for $-s_{k+1}<\eta<s_{k+1}$ with the  notation
$$a=a_{0k},\ b=a_{1k},\ c=a_{2k},$$
the  estimate
\begin{eqnarray*}
\big|Q_{k}(\eta)\big|&=&|a+b\eta+c{\eta^2}| = \big|\Psi^{(2)}_{k+1}(\xi_{0},\eta)-\eta-Q^{(2)}_{k}(\eta)\big|\\[0.2cm]
&\leq &\big|\Psi_{k+1}(\xi_{0},\eta)-\Omega_{k+1}(\xi_{0},\eta)-Q_{k}(\eta)\big|\\[0.2cm]
 & \leq &\big|\Psi_{k+1}-\Omega_{k+1}-Q_{k}\big|_{\mathbb{R}^2\cap D_{k+1}}\leq N
\end{eqnarray*}
holds by virtue of (3.25), where we  put
$$N=\Big({{5\over 24}+{{c_{1}}\over {18c_{1}+6c_{2}}}}\Big)M_{k+1}.$$
Then we have
$$|a+b\eta+c{\eta^2}|\leq N.$$
For $\eta=0$, we obtain $|a|\leq N$.\ Therefore we get
$$|b\eta+c{\eta^2}|\leq |a|+|a+b\eta+c{\eta^2}|\leq 2N.$$
Let $\eta=\pm \tilde{\sigma},\ 0<\tilde{\sigma}<s_{k+1}$,\ we get
$$|b|\tilde{\sigma}+|c| \tilde{\sigma}^2=|\pm b \tilde{\sigma}+c \tilde{\sigma}^2|\leq 2N.$$
Letting  $\tilde{\sigma}\rightarrow s_{k+1}$,\ we have
\begin{eqnarray*}
|a+b\eta+c{\eta^2}|&\leq& |a|+|b||\eta|+|c||{\eta^2}|\\[0.2cm]
&\leq&|a|+|b|s_{k+1}+|c|s^2_{k+1}
\leq 3N
\end{eqnarray*}
for all $\eta\in \mathbb{C},\ |\eta|<s_{k+1}$,\ hence
$$\big|Q_{k}\big|_{D_{k+1}}\leq 3N.$$
In the previous setting, we have obtained
$$\big|H_{k+1}-\Omega_{k+1}-Q_{k}\big|_{D_{k+1}}\leq \Big({{5\over 24}+{{c_{2}}\over {18c_{1}+6c_{2}}}}\Big)M_{k+1},$$
hence
\begin{eqnarray*}
\big|H_{k+1}-\Omega_{k+1}\big|_{D_{k+1}}&\leq& 3N+\Big({{5\over 24}+{{c_{2}}\over {18c_{1}+6c_{2}}}}\Big)M_{k+1}\\[0.2cm]
&\leq& {5\over 6}M_{k+1}+{{3c_{1}+c_{2}}\over {{18c_{1}+6c_{2}}}}M_{k+1}\\[0.2cm]
&\leq& {5\over 6}M_{k+1}+{1\over 6}M_{k+1}=M_{k+1}.
\end{eqnarray*}
This inequality obviously gives the wanted estimate $(3.10)_{k+1}$.

The proof by induction  for  justifying  the  iteration  process has finished.  It  remains  to  find  a better form  for  condition  (3.22).  Equivalently we may write
$$f(q)\ :=p \log(1+q)-\log(1-q)\leq(p-2\tau-1)\log2$$
for $0<q\leq{1\over 4},\ \ p>2\tau+1>2n+1>3$, we have ${d^{2}f\over {d{q^2}}}(q)\leq 0,$ hence
$$f(q)=f(q)-f(0)<q {df\over dq}(0)=q(p+1).$$
As a consequence
$$0<q\leq{{p-2\tau-1}\over {p+1}}\log2$$
is  sufficient  for  (3.22).  This is  one of the  conditions for  $q$ appearing  in  (\ref{b7}).

\section{Linear difference equations}
In  this  section we will solve the difference  equations
\begin{equation}\label{d1}
\begin{array}{ll}
u(x + \alpha, y)-u(x, y) = \varepsilon v(x, y) + f(x,y), \\[0.2cm]
v(x + \alpha, y)-v(x, y) = g(x, y)-[g](y),
\end{array}
\end{equation}
which plays a central  role  in  the  proof of the  inductive  theorem. Here the mean value of the function $g(x,y)$ over the  variable $x$ is defined  by $[g](y)=\lim \limits_{T\rightarrow\infty}{1 \over T}\int _{0}^{T}g(x,y)dx$, and $\alpha$ is  a real number satisfying  the  Diophantine inequalities
\begin{equation}\label{d2}
\begin{array}{ll}
\left\{\begin{array}{ll}
\Big|{\langle k,\omega \rangle {\alpha \over {2\pi}}-j}\Big|\geq {\gamma \over {|k|^\tau}},\ \ \ \ \mbox{for all}\ \  k \in \mathbb{Z}^n\backslash\{0\},\ \  j \in \mathbb{Z},\\[0.4cm]
 0<\gamma<{1\over 2},\ \tau>n.
 \end{array}\right.
\end{array}
\end{equation}
The functions $f(\cdot,y),g(\cdot,y)\in Q_{r}(\omega)$  are given holomorphic functions  of the  complex variables
$x,y $, and    $u ,  v$  are   wanted holomorphic functions  of the  complex variables
$x,y.$  $\varepsilon$ is  a positive constant  to  be determined  in  such a way
that  the  functions  $u,\ v$ will be of the  same size.

In  order to  get estimates  for  $u,v$ which are good enough for  the  proof of Theorem  \ref{thm2.6},
some technical  preparations  have to  be made.

\begin{lemma}[Lemma\ 3.3\ in {\cite{Russmann74}}]\label{lem4.1}
Let $\bar{\omega}=(\bar{\omega}_1,\bar{\omega}_2,\cdots,\bar{\omega}_\ell)\in{\mathbb{R}^\ell}$ satisfying the inequalities \ $D(k,\bar{\omega})\geq \psi(|k|)$, where \ $D(k,\bar{\omega})=\min \limits_{{j\in\mathbb{Z}}}\Big|\Big\langle(k,j),\bar{\omega}\Big\rangle\Big|, k\in {\mathbb{Z}^{\ell-1}\backslash \{0\}}$, $\psi$\ is an approximation function.\ Then for \ $m=1,2,\cdots,$\ we have
$$\sum \limits _{\substack{\bar{k}\in \mathbb{Z}^\ell\\ 0<|\bar{k}|\leq m}}{1\over {|\langle\bar{ k},\bar{\omega}\rangle|}^2}\leq {{\pi^2}\over 8}{{3^{\ell+2}\over {\psi(m)^2}}},$$
where $\bar{k}=(k,j)$.

\end{lemma}

If we choose
$$
\begin{array}{ll}
\bar{\omega}=(\omega_1{\alpha\over{2\pi}},\omega_2{\alpha\over{2\pi}},\cdots,\omega_n{\alpha\over{2\pi}},\ -1),\ \ell=n+1,\ \psi(t)=\gamma t^{-\tau}, \\[0.2cm]
k=(k_1,k_2,\cdots,k_n)\in \mathbb{Z}^n\backslash\{0\},\ j\in \mathbb{Z},\ \bar{k}=(k,j),
\end{array}
$$
then by Lemma \ref{lem4.1} and the Diophantine inequalities (\ref{d2}), we obtain
$$\sum \limits _{\substack{k\in \mathbb{Z}^n\\ 0<|k|\leq m}}{1\over {\Big|{\langle k,\omega \rangle {\alpha \over {2\pi}}-j}\Big|}^2}\leq {{\pi^2}\over 8}{3^{n+3}}{\gamma^{-2}}m^{2\tau}.$$
Meanwhile,
$$
\Big|e^{i\langle k,\omega\rangle \alpha}-1\Big|\ge \pi\Big|{\langle k,\omega \rangle {\alpha \over {2\pi}}-j}\Big|.
$$
Therefore
\begin{equation}\label{d3}
\sum \limits _{\substack{k\in \mathbb{Z}^n\\ 0<|k|\leq m}}{1\over{\Big|e^{i\langle k,\omega\rangle \alpha}-1\Big|}^2} \leq {{3^{n+3}}\over 8}{\gamma^{-2}} m^{2\tau}.
\end{equation}

\begin{lemma}\label{lem4.2}
For $r  > 0$, let $f : \big\{x\in \mathbb{C}: |\Im\ x|< r\big\}\mapsto \ \mathbb{C}$
be a holomorphic function  and $f\in Q_{r}(\omega)$.  Then we have the  estimate
$$
\sum \limits _{k\in \mathbb{Z}^n}|f_{k}|^2{e^{2r|k|}}\leq {2^n}\big|f\big|_{r}^2,
$$
where
$$f_{k}={1 \over {(2\pi)}^n}\int _{\mathbb{T}^n}F(\theta)e^{-i\langle k,\theta\rangle}d\theta,\ \  k\in \mathbb{Z}^n$$
are the  Fourier coefficients  of $f$,\  $F(\theta)$ is the shell function of $f$ and
$$\big|f\big|_{r}=\sup \limits_{\theta\in \Pi_{r}^{n}}\Big|\sum_{k}f_{k}e^{i\langle k,\theta\rangle}\Big|=\sup \limits_{\theta\in \Pi_{r}^{n}}|F(\theta)|.$$
\end{lemma}
\noindent\textbf{Proof}:
The Fourier coefficients $f_{k}$ of $f$ are given by
$$f_{k}={1 \over {(2\pi)}^n}\int _{\mathbb{T}^n}F(\theta)e^{-i\langle k,\theta\rangle}d\theta,\ \ k=(k_1,k_2,\cdots,k_n)\in \mathbb{Z}^n.$$
For every $\lambda\in \mathbb{R}^n$ with $|\lambda|=\max \limits_{1\leq j\leq n} |\lambda_{j}|<r,$  the function $x\mapsto f(x+i\lambda)$ which domain is $r-|\lambda|,\ $ its Fourier coefficients are
$$f_{k}(\lambda)={1 \over {(2\pi)}^n}\int _{\mathbb{T}^n}F(\theta+i\lambda)e^{-i\langle k,\theta\rangle}d\theta,\ \ k\in \mathbb{Z}^n.$$
By Bessel's inequality,
$$\sum \limits _{k\in \mathbb{Z}^n}|f_{k}(\lambda)|^2\leq {1 \over{(2\pi)}^n}\int _{\mathbb{T}^n}|F(\theta+i\lambda)|^2d\theta.$$
Hence
\begin{equation}\label{d5}
\sum \limits _{k\in \mathbb{Z}^n}|f_{k}(\lambda)|^2\leq \big|f\big|_{r}^2,\ \ \ \ |\lambda|<r.
\end{equation}

Define a new function
$$\lambda\mapsto f_{k}(\lambda)e^{\langle k,\lambda\rangle}={1 \over {(2\pi)}^n}\int _{\mathbb{T}^n}F(\theta+i\lambda)e^{-i\langle k,\theta+i\lambda\rangle}d\theta,$$
then
$${\partial\over \partial \lambda_{j}}\Big({f_{k}(\lambda)e^{\langle k,\lambda\rangle}}\Big)={1 \over {(2\pi)}^n}\int _{\mathbb{T}^n}i{\partial\over \partial \theta_{j}}F(\theta+i\lambda)e^{-i\langle k,\theta+i\lambda\rangle}d\theta.$$
Since $F(\theta+i\lambda)e^{-i\langle k,\theta+i\lambda\rangle}$ is $2\pi$-periodic in  $\theta_{j}\ (j=1,2,\cdots,n)$,\ then
$${\partial\over \partial \lambda_{j}}\Big({f_{k}(\lambda)e^{\langle k,\lambda\rangle}}\Big)=0.$$
Hence the function ${f_{k}(\lambda)e^{\langle k,\lambda\rangle}}$ is independent of $\lambda$, and
$${f_{k}(\lambda)e^{\langle k,\lambda\rangle}}=f_{k}(0)=f_{k},$$
and consequently
$${|f_{k}(\lambda)|^2}e^{2{\langle k,\lambda\rangle}}=|f_{k}|^2.$$
Finally,\ by (\ref{d5}), we have
\begin{equation}\label{d6}
\sum \limits _{k\in \mathbb{Z}^n}|f_{k}|^2e^{-2{\langle k,\lambda\rangle}}=\sum \limits _{k\in \mathbb{Z}^n}|f_{k}(\lambda)|^2\leq \big|f\big|_{r}^2,\ \ \ \ |\lambda|<r.
\end{equation}

Define $e_i\in \mathbb{R}^n\ (i=1,2,\cdots,2^n)$ which have components $\pm 1$, and
$$\mathbb{Z}_{i}=\Big\{k\in \mathbb{Z}^n : \langle k, e_i\rangle=-|k|\Big\}.$$
Then
\begin{equation}\label{d7}
\bigcup_{i=1}^{2^n}\mathbb{Z}_{i}=\mathbb{Z}^n.
\end{equation}
Let $\lambda=s e_i$ in (\ref{d6}), we obtain
$$\sum \limits _{k\in \mathbb{Z}_{i}}|f_{k}|^2 e^{2s|k|}\leq  \big|f\big|_{r}^2,\ \ \ \ 0<s<r,\ \  i=1,2,\cdots,2^n.$$
Passing to  the limit $s\rightarrow r$ yields
$$\sum \limits _{k\in \mathbb{Z}_{i}}|f_{k}|^2 e^{2r|k|}\leq  \big|f\big|_{r}^2.$$
Adding these inequalities and by (\ref{d7}),\ we have
$$\sum \limits _{k\in \mathbb{Z}^n}|f_{k}|^2{e^{2r|k|}}\leq {2^n}\big|f\big|_{r}^2.$$
The proof of this lemma is completed. \qed

\begin{definition}
(i) For $D=D(r,s)$, denote by $P^{m}(r,s)=P^{m}(D)\ (m=1,2,\cdots)$ the  linear  space of all holomorphic functions $f:D \rightarrow \mathbb{C}^{m}$ satisfying
$$f\circ \sigma \big|_ D=f,\ \ f(\mathbb{R}^2)\subseteq \mathbb{R}^m.$$
Clearly $f,g \in T(D)$ implies $f-g \in P^2(D).$\\
(ii) For a function $f\in P^{m}(D)$, denote its  mean value  over the  variable $x$   by
$$[f](y)=\lim \limits_{T\rightarrow\infty}{1 \over T}\int _{0}^{T}f(x,y)dx.$$
\end{definition}

\begin{theorem}\label{thm4.4}
Let $\alpha$ be a real  number satisfying (\ref{d2}),\ and $f(\cdot,y)\in Q_{r}(\omega)$  be a function
belonging to  $P^{1}(r,s)$  for  some positive constants  $r,s.$  Then the  difference equation
$$
u(x + \alpha, y)-u(x, y) = f(x, y) - [f](y)
$$
 has a unique solution $u\in Q(\omega),\ u\in P^{1}(r,s)$  with $[u]=0.$  For this  solution  the  estimate
 \begin{equation}\label{d9}
\big|u\big|_{r-\rho,s}\leq \varepsilon^{-1}\big|f\big|_{r,s},\ \ \ \ 0<\rho<r
\end{equation}
holds, where $\varepsilon$ is  defined by
\begin{equation}\label{d10}
\varepsilon=\varepsilon(\rho)=6^{-{{n+1}\over 2}}\  {\gamma\over {\Gamma(\tau+1)}} \ \rho^{\tau}.
\end{equation}

\end{theorem}
\noindent\textbf{Proof}:
Since the  restriction  of $f(x,y)$ onto $\mathbb{R}^2$ is  a continuously differentiable  and quasi-periodic function in $x$, it  can be expanded into  its  Fourier series
$$
f(x,y)=\sum \limits _{k\in \mathbb{Z}^n} f_{k}(y)e^{i \langle k,\omega \rangle x},
$$
where
$$
f_{k}(y)={1 \over {(2\pi)}^n}\int _{\mathbb{T}^n}F(\theta,y)e^{-i\langle k,\theta\rangle}d\theta,\ \ k=(k_1,k_2,\cdots,k_n)\in \mathbb{Z}^n\\[0.2cm]
$$
are the  Fourier coefficients  of $f(x,y)$.\, An application of Lemma \ref{lem4.2} to  the restriction  of $f(x,y)$  onto $D(r,s)$
yields that
\begin{equation}\label{d13}
\sum \limits _{k\in \mathbb{Z}^n}|f_{k}(y)|^2{e^{2r|k|}}\leq {2^n}\big|f\big|_{r,s}^2.
\end{equation}

Let
$$u(x,y)=\sum \limits _{k} u_{k}(y)e^{i \langle k,\omega \rangle x}.$$
After straightforward calculations we obtain the relation between Fourier coefficients $f_{k}(y)$ and $u_{k}(y)$ as follows
$$
u_{k}(y)={{f_{k}(y)}\over {e^{i\langle k, \omega\rangle \alpha}-1}},\ \ \ k\neq 0,
$$
then
$$
u(x,y)=\sum \limits _{k\neq 0} {{f_{k}(y)}\over {e^{i\langle k, \omega\rangle \alpha}-1}}e^{i \langle k,\omega \rangle x},
$$
which is the  uniquely determined Fourier expansion of the  wanted solution $u$  satisfying $u \in Q(\omega)$ with
$[u]=0.$

Firstly,\ we estimate the sum
$$g_{m}(y)=\sum \limits _{1\leq |k|\leq m} \Big|{{f_{k}(y)}\over {e^{i\langle k, \omega\rangle \alpha}-1}}\Big|e^{|k|r}.$$
By Cauchy-Schwarz inequality,\ we have
$$g_{m}(y)\leq \sqrt{\sum \limits _{k\in \mathbb{Z}^n}{\big|f_{k}(y)\big|^2}e^{2|k|r}}\sqrt{ \sum \limits _{\substack{k\in \mathbb{Z}^n\\ 0<|k|\leq m}}{\Big|e^{i\langle k,\omega\rangle \alpha}-1\Big|}^{-2}}.$$
By\ $(\ref{d3}),\ (\ref{d13}),$ we obtain
$$\sqrt{\sum \limits _{k\in \mathbb{Z}^n}{\big|f_{k}(y)\big|^2}e^{2|k|r}}\leq 2^{n\over 2}\  \big|f\big|_{r,s},$$
$$\sqrt{ \sum \limits _{\substack{k\in \mathbb{Z}^n\\ 0<|k|\leq m}}{\Big|e^{i\langle k,\omega\rangle \alpha}-1\Big|}^{-2}}\leq {3^{{n\over 2}+{3\over 2}}\over {2\sqrt 2}}\ {{m^\tau}\over \gamma},$$
hence
$$g_{m}(y)\leq {{3\sqrt 3}\over {2\sqrt 2}}\ 6^{n\over 2}\ {{m^\tau}\over \gamma}\ \big|f\big|_{r,s}\leq 6^{{n+1}\over 2}\ {{m^\tau}\over \gamma}\ \big|f\big|_{r,s}.$$
Set $g_{0}(y)=0$, we get
$$\sum \limits _{0< |k|\leq N} \Big|{{f_{k}(y)}\over {e^{i\langle k, \omega\rangle \alpha}-1}}\Big|e^{|k|(r-\rho)}
=(1-e^{-\rho})\sum \limits _{m=1}^{N}g_{m}(y)e^{-m\rho}+g_{N}(y)e^{-(N+1)\rho}.$$
Letting $N\rightarrow \infty $,\ we have
\begin{eqnarray*}
\sum \limits _{k\neq 0} \Big|{{f_{k}(y)}\over {e^{i\langle k, \omega\rangle \alpha}-1}}\Big|e^{|k|(r-\rho)}
&\leq & (1-e^{-\rho})\sum \limits _{m=1}^{\infty}g_{m}(y)e^{-m\rho}\\[0.2cm]
&\leq & {6^{{n+1}\over 2}\over \gamma}\big|f\big|_{r,s}\sum \limits _{m=1}^{\infty}{m^\tau}(e^{-m\rho}-e^{-(m+1)\rho}).
\end{eqnarray*}
For the last series we get the estimate
\begin{eqnarray*}
\sum \limits _{m=1}^{\infty}{m^\tau}\Big(e^{-m\rho}-e^{-(m+1)\rho}\Big)
&=&\rho\sum \limits _{m=1}^{\infty}{m^\tau}\int_{m}^{m+1}e^{-x\rho}d x\\[0.2cm]
&\leq & \rho \sum \limits _{m=1}^{\infty} \int_{m}^{m+1} x^{\tau}e^{-x\rho}d x=\rho \int_{1}^{+\infty}x^{\tau}e^{-x\rho}d x\\[0.2cm]
&\leq & \rho^{-\tau} \int_{1}^{+\infty}(x \rho)^{\tau}e^{-x\rho}d (x\rho)\leq \rho^{-\tau} \int_{0}^{+\infty}t^{\tau}e^{-t}d t\\[0.2cm]
&= & \rho^{-\tau} \Gamma(\tau+1).
\end{eqnarray*}
Hence
$$\big|u\big|_{r-\rho,s}\leq 6^{{n+1}\over 2}\ {\Gamma(\tau+1)\over \gamma}\ \rho^{-\tau}\ \big|f\big|_{r,s},$$
which completes the proof of the lemma.\qed

Now we are ready  to  solve  equation (\ref{d1}).

\begin{theorem}\label{thm4.5}
Let $\alpha$ be a real  number satisfying (\ref{d2}),\ and $f(\cdot,y),g(\cdot,y)\in Q_{r}(\omega)$  be  functions
belonging to  $P^{1}(r,s)$  and satisfying  the  estimates
\begin{equation}\label{d16}
\big|f\big|_{r,s}\leq M,\ \ \ \ \big|g\big|_{r,s}\leq M
\end{equation}
with some positive constants  $r,\ s,\ M.$ Then the  difference equations (\ref{d1}) with $\varepsilon$ defined in (\ref{d10}) have a unique solution $u,v \in Q(\omega),$ $u,v\in P^{1}(r,s)$ with $[u]=0.$ For this  solution the  estimates
\begin{equation}\label{d17}
\big|u\big|_{r-2\rho,s}\leq 2\varepsilon^{-1}M,
\end{equation}

\begin{equation}\label{d18}
\big|v\big|_{r-\rho,s}\leq 2\varepsilon^{-1}M
\end{equation}
 are valid for $0<2\rho<r$.
\end{theorem}
\noindent\textbf{Proof}:
 In   the first equation of (\ref{d1}) the  mean value  must vanish on both sides. Hence we get the condition
\begin{equation}\label{d19}
[v]= -\varepsilon^{-1}[f]
\end{equation}
for  the  mean value of $v$. As a consequence,  we have $[v]\in P^{1}(r,s)$ and
\begin{equation}\label{d20}
\big|[v]\big|_{r,s}\leq \varepsilon^{-1}M
\end{equation}
in  view of (\ref{d16}). Theorem \ref{thm4.4} gives a unique solution $v=\tilde{v}\in P^{1}(r,s)$ of the second equation of (\ref{d1}) with $[\tilde{v}]=0.$ This solution  has the  estimate
\begin{equation}\label{d21}
\big|\tilde{v}\big|_{r-\rho,s}\leq \varepsilon^{-1}M
\end{equation}
because of (\ref{d16}). Define $v=\tilde{v}+[v]$, we obtain the  uniquely determined solution $v\in P^{1}(r,s)$ of the second equation of (\ref{d1})  satisfying
(\ref{d19}).\ This solution  has the  estimate (\ref{d18}) as a consequence  of (\ref{d20}) and (\ref{d21}).

Define $h=\varepsilon \tilde{v}+f$, note that  $\tilde{v}$ is defined in $D(r-\rho,s)$,\ then  $h$ is well defined in $D(r-\rho,s)$. As a consequence we have
\begin{equation}\label{d22}
\big|h\big|_{r-\rho,s}=\big|\varepsilon\tilde{v}+f\big|_{r-\rho,s}\leq 2M
\end{equation}
and
\begin{eqnarray*}
h(x,y)-[h](y)
&=&\varepsilon\tilde{v}+f-[\varepsilon \tilde{v}+f] \\[0.2cm]
&=&  \varepsilon\tilde{v}+f-\varepsilon [\tilde{v}]-[f] \\[0.2cm]
&=&  \varepsilon\tilde{v}+f+\varepsilon[v]\\[0.2cm]
&= & \varepsilon v+f.
\end{eqnarray*}
Hence, the first equation of (\ref{d1}) can  be rewritten in  the  form

\begin{equation}\label{d23}
u(x + \alpha, y)-u(x, y) = h(x, y) - [h](y).
\end{equation}
Thus Theorem  \ref{thm4.4} gives a uniquely determined solution  $u\in P^{1}(r,s)$ of (\ref{d23}) with $[u]=0.$ For an estimate of $u$ we apply Theorem  \ref{thm4.4} to  (\ref{d23}) restricted  to $D(r-\rho,s)$ such that  in  (\ref{d9}) we have to  replace $f$ by $h$ and $r$ by $r-\rho.$
Then (\ref{d17}) follows by means of (\ref{d22}). The proof is finished. \qed

\section{The inductive  theorem }
First  of all  we give together  constants,  domains, and mappings appearing in  the
formulation  of  the  inductive  theorem.\\
\vspace{0.2cm}
(I)  Constants and their  relations
\vspace{0.2cm}

We introduce  the  constants
\begin{equation}\label{e1}
\omega,\ \gamma,\ \tau,\ M,\ q,\ \varepsilon,\ \varepsilon_{+},\ r,\ r_{+},\ s,\ s_{+},\ r', r'_{+},\ s,\ s'_{+}
\end{equation}
and the  auxiliary  constants  $\theta,\rho$ satisfying  the  relations
$$
\begin{array}{ll}
\Big|{\langle k,\omega \rangle {\alpha \over {2\pi}}-j}\Big|\geq {\gamma \over {|k|^\tau}},\ \ \ \mbox{for all}\ \  k \in \mathbb{Z}^n\backslash\{0\},\ \  j \in \mathbb{Z},\\[0.2cm]
0<\gamma<{1\over 2},\ \ \ n\leq \tau,\ \ \ 0<r\leq 1,\\[0.2cm]
0<q\leq\big({\theta\over 10}\big)^2,\ \ \  s={\theta\rho \over 50},\ \ \  \theta=2^{-\tau},\ \ \ \  \rho={r\over 6},\\[0.2cm]
\varepsilon=6^{-{{n+1}\over 2}}\  {\gamma\over {\Gamma(\tau+1)}} \ \rho^{\tau},\ \ \ M={1\over 3}q\varepsilon s,\ \ \ \\[0.2cm]
 r'={4\over 3}(r-s),\ \ \ s'={4\over 3}s,\\[0.2cm]
{r_+\over r}={s_+\over s}={r'_+\over r'}={s'_+\over s'}={1\over 2},\ \ \ \ {\varepsilon_+\over \varepsilon}=\theta.
\end{array}
$$

\vspace{0.2cm}
(II)  Domains and Mappings
\vspace{0.2cm}

Choose
$$D=D(r,s),\ D_{+}=D(r_{+},s_{+}),\ D'=D(r',s'),\ D'_{+}=D(r'_{+},s'_{+}),$$
$$D_{+}^{\star}=D(r'_{+}-{1\over 7}s,s'_{+}-{1\over 7}s)=D'_{+}-{1\over 7}s,$$
and introduce  the  mappings
\begin{equation}\label{e2}
\Omega(x,y):\quad
x_1=x+\alpha+{\varepsilon}y,\quad
y_1=y.
\end{equation}
\begin{equation}\label{e3}
{\Omega}_{+}(x,y):\quad
x_1=x+\alpha+{\varepsilon_{+}}y,\quad
y_1=y.
\end{equation}
\begin{equation}\label{e4}
 \Theta(x,y):\quad
x_1=x,\quad
y_1=\theta y
\end{equation}
for  all $(x,y)\in \mathbb{C}^2,$ where we use the  same symbol  for  the  mappings $\Omega,\ \Omega_{+}, \Theta$ as well
as for  their  restrictions  to  subsets  of $\mathbb{C}^2.$

For the proof of the inductive theorem we need two useful lemmas.
\begin{lemma} [Lemma\ 5\ in {\cite{Russmann83}}]\label{lem5.1}
Let $D\subseteq \mathbb{C}^\ell$ be an open and convex set,\ and $F:D\mapsto \mathbb{C}^m$ be a holomorphic function. Then for any $d>0$, we have
$$|F(x)-F(y)|\leq {1\over d}|x-y|\sup \limits_{x\in D}|F(x)|,\ \ \  \mbox{for all}\ \  x,y\in D-d.$$
\end{lemma}

\begin{lemma} [Lemma\ 6\ in {\cite{Russmann83}}]\label{lem5.2}
For some $r>0,$ let $f:\big\{z\in \mathbb{C}:|z|<r\big\}\mapsto \mathbb{C}$ be a holomorphic function with power series expansion
$$f(z)=\sum_{k=0}^{\infty}f_{k}z^k.$$
Then for the polynomial
$$f_{m-1,q}(z)=\sum_{k=0}^{m-1}\Big(1-q^{2(m-k)}\Big)f_{k}z^k,\ \ \ \ m=1,2,\cdots$$
of degree $m-1$ depending on $q \ ( 0<q<1)$, we have the estimate
$$\sup_{|z|\leq q r}\big|f(z)-f_{m-1,q}(z)\big|\leq q^m\sup_{|z|<r}|f(z)|.$$
\end{lemma}

\begin{theorem} [Inductive Theorem]\label{thm5.3}
Let constants (\ref{e1}) and auxiliary constants $\theta,\rho$ be given such that  the  relations in  (I)  are satisfied,  and let  domains $D,\ D_{+},\ D',\ D'_{+},\ D_{+}^{\star}$ and mappings $\Omega,\ \Omega_{+},\ \Theta$ be given as in  (II).  Then for  any mapping
$$H\ : \ D\rightarrow D',\ \ \ \ H\in T(D), $$
$H-\Omega$ is quasi-periodic with the frequency $\omega$ in the first variable and satisfying
\begin{equation}\label{e5}
\big|H-\Omega\big|_{D}\leq M,
\end{equation}
there  are mappings $W\in T(D'_{+})$, $\Phi_{+}\in T(D_{+})$ such that  the  diagram
$$\xymatrix{
D_{+} \ar [rr]^{W\big |_{D_{+}}} \ar "2,2"_{\Phi_{+}}&{}&
  D\ar [dd]^{H}\\
{}&{D^\star_{+}}\ar "3,1"_{id}&{}\\
 {D'_{+} }\ar [rr]_{W}&{}
&{D'}
}$$
exists and commutes.  Moreover $W-\Theta, \Phi_{+}-\Omega_+-Q$ are quasi-periodic with the frequency $\omega$ in the first variable and the  following  estimates are satisfied
\begin{equation}\label{e6}
\big|W-\Theta \big|_{D'_{+}}\leq {2\over 3}qs,
\end{equation}
\begin{equation}\label{e7}
\theta(1-q)|\zeta-\zeta'|\leq |W(\zeta)-W(\zeta')|\leq (1+q)|\zeta-\zeta'|,\ \  \zeta,\ \zeta' \in D'_{+},
\end{equation}
\begin{equation}\label{e8}
\big|\Phi_{+}-\Omega_{+}-Q\big|_{D_{+}}\leq {5\over 48}\theta M,
\end{equation}
where $Q : D_{+} \mapsto \mathbb{C}^2$ is defined  by
$$(\xi,\eta)\mapsto Q(\eta)=(0,a_{0}+a_{1}\eta+a_{2}\eta^2)$$
with some constants $a_{0},\ a_{1},\ a_{2}\in \mathbb{R}.$
\end{theorem}
\noindent\textbf{Proof}:\, 1) First  of all, define
$$h=\Big(\begin{matrix} f \\ g \end{matrix} \Big)\ :=H-\Omega,$$
and by assumption, $H-\Omega$ is quasi-periodic with the frequency $\omega$ in the first variable and $h\in P^2(D).$

The results  of Section  4 enable us to  solve the  linear  difference  equations
\begin{equation}\label{e9}
\left\{\begin{array}{ll}
u(x + \alpha, y)-u(x, y) = \varepsilon v(x, y) + f(x,\theta y),\\[0.4cm]
v(x + \alpha, y)-v(x, y) = g(x, \theta y)-[g](\theta y).
 \end{array}\right.
\end{equation}

Let $d\, \Omega$ be the differential of $\Omega$ and define
$$w=\Big(\begin{matrix} u \\ v \end{matrix} \Big),$$
\begin{equation*}\label{e10}
{\Omega}^{\star}(x,y):\quad
x_1=x+\alpha,\quad
y_1=y,
\end{equation*}

\begin{equation*}\label{e11}
h^{\star}(x,y):\quad
x_1=0,\quad
y_1=[g](y).
\end{equation*}
With these definitions and notations, we have
$$w\circ \Omega^{\star}=\Bigg(\begin{matrix} u(x + \alpha, y) \\ v(x + \alpha, y) \end{matrix} \Bigg),$$
$$(d\, \Omega) w=\Bigg(\begin{matrix} 1\ \ \ \varepsilon \\ 0\ \ \ 1 \end{matrix} \Bigg)\Bigg(\begin{matrix} u(x, y) \\ v(x, y) \end{matrix} \Bigg)=\Bigg(\begin{matrix} u(x, y)+\varepsilon v(x, y) \\ v(x, y) \end{matrix} \Bigg),$$
$$h\circ \Theta=\Bigg(\begin{matrix} f \\ g \end{matrix} \Bigg)\circ \Theta=\Bigg(\begin{matrix} f(x,\theta y) \\ g(x,\theta y) \end{matrix} \Bigg),$$
$$h^{\star}\circ \Theta=\Bigg(\begin{matrix} 0 \\ [g] \end{matrix} \Bigg)\circ \Theta=\Bigg(\begin{matrix} 0 \\ [g](\theta y) \end{matrix} \Bigg).$$
Hence,  the difference  equations (\ref{e9}) can  be written  in  the  more compact form
\begin{equation}\label{e12}
w\circ \Omega^{\star}= (d\, \Omega) w+ h \circ \Theta -h^{\star} \circ \Theta.
\end{equation}

Define
$$w=\Big(\begin{matrix} u \\ v \end{matrix} \Big)\ :=W-\Theta,\ \ \ \phi\ :=\Phi_{+}-\Omega_{+},$$
we ought to  show $w\in P^2(D),\ \phi\in P^2(D_{+}).$

After having obtained a solution $w\in P^2(D'_{+})$ of (\ref{e12}), we try  to  determine $\phi$ from  the  equation
\begin{equation}\label{e13}
H \circ W\big |_{D_{+}}=W \circ \Phi_+,
\end{equation}
which holds because the  diagram in the inductive theorem  can commute.

Now (\ref{e13}) can be rewritten in the form
\begin{equation}\label{e14}
(h+\Omega)\circ (\Theta+w)=(\Theta+w)\circ (\phi+\Omega_{+}),
\end{equation}
which is
$$h\circ (\Theta+w)+\Omega\circ (\Theta+w)=w\circ (\phi+\Omega_{+})+\Theta\circ (\phi+\Omega_{+}).$$
First,
\begin{eqnarray*}
\Omega\circ (\Theta+w)
&=&\Bigg(\begin{matrix} x+u(x,y)+\alpha+\varepsilon(\theta y+v) \\ \theta y+v \end{matrix} \Bigg)\\[0.2cm]
&=&\Bigg(\begin{matrix} x+\alpha+\varepsilon\theta y \\ \theta y \end{matrix} \Bigg)+\Bigg(\begin{matrix} u+\alpha+\varepsilon v \\ v \end{matrix} \Bigg)-\Bigg(\begin{matrix} \alpha \\ 0 \end{matrix} \Bigg)\\[0.2cm]
&=& \Omega\circ \Theta+\Omega\circ w-\Bigg(\begin{matrix} \alpha \\ 0 \end{matrix} \Bigg),
\end{eqnarray*}
$$\Theta\circ (\phi+\Omega_{+})=\Theta\circ \phi+\Theta\circ \Omega_{+}.$$
Since ${{\varepsilon_{+}}\over {\varepsilon}}=\theta$,\ then
$$\Omega\circ \Theta=\Bigg(\begin{matrix} x+\alpha+\varepsilon\theta y \\ \theta y \end{matrix} \Bigg)=\Bigg(\begin{matrix} x+\alpha+\varepsilon_{+} y \\ \theta y \end{matrix} \Bigg)$$
and
$$\Omega\circ w-\Bigg(\begin{matrix} \alpha \\ 0 \end{matrix} \Bigg)=\Bigg(\begin{matrix} u+\varepsilon v \\ v \end{matrix} \Bigg)=(d\, \Omega) w,$$
$$\Theta\circ \Omega_{+}=\Theta\circ\Bigg(\begin{matrix} x+\alpha+\varepsilon_{+} y \\  y \end{matrix} \Bigg)=\Bigg(\begin{matrix} x+\alpha+\varepsilon_{+} y \\ \theta y \end{matrix} \Bigg),$$
which implies that
$$\Omega\circ \Theta=\Theta\circ \Omega_{+}.$$
Thus (\ref{e13}) is changed into the form
$$\Theta\circ\phi=h\circ (\Theta+w)-w\circ (\phi+\Omega_{+})+(d\, \Omega) w.$$
If we define
$$\phi={\Theta^{-1}}(z+h^{\star}\circ\Theta),$$
which leads to
$$\Theta\circ\phi=z+h^{\star}\circ\Theta,$$
then
$$z+h^{\star}\circ\Theta =h\circ (\Theta+w)-w\circ (\phi+\Omega_{+})+(d\, \Omega) w.$$

Now (\ref{e13}) gets by  (\ref{e12}) the  form
$$
\begin{array}{ll}
 \left \{\begin{array}{ll}
F_{1}=w\circ \Omega_{+}-w\circ(\Omega_{+}+\phi),\\[0.2cm]
F_{2}=w\circ \Omega^{\star}-w\circ \Omega_{+},\\[0.2cm]
F_{3}=h\circ (\Theta+w)-h\circ\Theta,\\[0.2cm]
z=F(z)\ :=F_{1}+F_{2}+F_{3}.
 \end{array}\right.
\end{array}
$$
Careful estimates will lead  to  a solution $z \in P^2(D'_{+})$ and $z$ is quasi-periodic with the frequency $\omega$ in the first variable , which implies that $\phi$ is quasi-periodic with the frequency $\omega$ in the first variable and $\phi \in P^2(D'_{+})$ can be  determined.

\medskip

2) Properties of  $W$.
\vspace{0.2cm}

Since $H-\Omega$ is quasi-periodic with the frequency $\omega$ in the first variable and $|H-\Omega|_{D}\leq M,$ using Theorem \ref{thm4.5} with $\rho={r\over 6}$, we get a solution $w$ which is quasi-periodic with the frequency $\omega$ in the first variable of (\ref{e12}) with
\begin{equation}\label{e15}
w\in P^{2}(r,t),\ \ \ \big|w\big|_{4\rho,t}\leq 2\varepsilon^{-1} M,\ \ \ t:={s\over \theta}.
\end{equation}
Since  $D'_{+}\subseteq D(4\rho,t)$, we define $W=\Theta+w|_{D'_{+}}$ and obtain $W-\Theta$ is quasi-periodic with the frequency $\omega$ in the first variable with
$$W\in T(D'_{+}),\ \ \ \big|W-\Theta\big|_{D'_{+}}\leq 2\varepsilon^{-1} M={2\over 3}qs,$$
which is the wanted estimate (\ref{e6}).

From (\ref{e15}) we also have
\begin{equation}\label{e16}
w\circ \Theta^{-1}\in P^{2}(r,s),\ \ \ \big|w\circ \Theta^{-1}\big|_{4\rho,s}\leq{2\over 3}qs.
\end{equation}
An application of Lemma \ref{lem5.1} with $d={2\over 3}s$ yields
$$\big|w\circ \Theta^{-1}(z)-w\circ \Theta^{-1}(z')\big|\leq d^{-1}\ |z-z'|\ \big|w\circ \Theta^{-1}\big|_{4\rho,s},\ \ \ \ z,z'\in D(r'_{+},{s\over 3}),$$
and which gets by means of (\ref{e16}) the form
$$\big|w\circ \Theta^{-1}(z)-w\circ \Theta^{-1}(z')\big|\leq q|z-z'|,\ \ \ \ z,z'\in D(r'_{+},{s\over 3}),$$
and consequently
$$(1-q)|z-z'|\leq \big|(\mbox{id}+w\circ \Theta^{-1})(z)-(\mbox{id}+w\circ \Theta^{-1})(z')\big|\leq (1+q)|z-z'|$$
for all $z,z'\in D(r'_{+},{s\over 3})$. Inserting $z=\theta \xi,\ z'=\theta\xi'$ for $\xi,\xi'\in D'_{+}$, we get (\ref{e7}) because of
$$\theta|\xi-\xi'|\leq |z-z'|\leq |\xi-\xi'|\ \ \text{and }\ \  \Theta(\xi')\in D'_{+}\subseteq D(r'_{+},{s\over 3}).$$

Now we look for the range of $W$. Since
$$qc\leq(1+q)c\leq {3\over 4}s\leq{r\over 4}\leq{r'\over 2},\ \ \ c:={2\over 3}s,$$
then
$$W(D'_{+})\subseteq D(r'_{+}+qc,\theta s'_{+}+qc)\subseteq D({r'\over 2}+{r'\over 2},{s'\over 4}+{3\over 4}s)\subseteq D',$$
$$W(D_{+})\subseteq D(r_{+}+qc,\theta s_{+}+qc)\subseteq D({r\over 2}+{r\over 2},{s\over 4}+{3\over 4}s)-c=D-{2\over 3}s\subseteq D,$$
$$\Theta(D_{+})= D(r_{+},\theta s_{+})\subseteq D(r_{+}+qc,\theta s_{+}+qc)\subseteq D-{2\over 3}s\subseteq D.$$

\medskip

3) Estimate for $F_{2}=w\circ \Omega^{\star}-w\circ \Omega_{+}$.

\vspace{0.2cm}

First $\varepsilon<1$ because $\Gamma(\tau+1)\geq 1$ for $\tau>n$. Hence we get
$$\Omega_{+}(D_{+})\subseteq D(r_{+}+\theta \varepsilon s_{+},s_{+})\subseteq D(r_{+}+{\theta \over 2}s,s_{+}).$$
Since $\Omega^{\star}(D_{+})=D_{+},\ s_{+}<t=\theta^{-1}s,$ we have
\begin{equation}\label{e17}
\Omega^{\star}(D_{+}),\ \Omega_{+}(D_{+})\subseteq D(4\rho-R,t),\ R=\rho-{\theta\over 2}s.
\end{equation}
Moreover applying Lemma \ref{lem5.1} to  $\xi\mapsto w(\xi,\eta)$ with $d=R,\ D=\big\{\xi\in\mathbb{C} : |\Im\ \xi|<4\rho\big\}$,\ we obtain
$$|w(\xi,\eta)-w(\xi',\eta)|\leq 2\varepsilon^{-1}M R^{-1}|\xi-\xi'|,\ \ \ (\xi,\eta),\ (\xi',\eta)\in D(4\rho-R,t).$$
Using (\ref{e17}) and the definitions of $\Omega^{\star},\ \Omega_{+}$ yields
$$\big|F_{2}\big|_{D_{+}}\leq 2\varepsilon^{-1}M R^{-1}\varepsilon_{+} s_{+}={s\over R}\theta M.$$
By (\ref{e17}) and (I) we have
$$R=\rho-{\theta\over 2}s={50\over \theta}s-{\theta s\over 2}> 49{s\over \theta},\ \ \ {{\theta s}\over 2}<{s\over \theta},$$
hence
\begin{equation}\label{e18}
F_{2}\in P^{2}(D_{+}),\ \ \ \big|F_{2}\big|_{D_{+}}\leq {1\over 49}\theta^2 M.
\end{equation}

\medskip
4) Estimate for $F_{3}=h\circ W\big|_{D_{+}}-h\circ\Theta$.
\vspace{0.2cm}

Since $\big|h\big|\leq M$, we apply Lemma \ref{lem5.1} with $d={2\over 3}s$  to obtain
$$|h(z)-h(z')|\leq {{3M}\over {2s}}|z-z'|,\ \ \ z,z'\in D-{2\over 3}s.$$
In the previous setting we have  proved $W(D_{+})\subseteq D-{2\over 3}s\subseteq D,\ \Theta(D_{+})\subseteq D-{2\over 3}s\subseteq D,\ $ then
$$\big|F_{3}\big|_{{D_{+}}}\leq {{3M}\over {2s}}\big|W-\Theta\big|_{{D_{+}}},\ \ \ \ M={1\over 3}q\varepsilon s.$$
Moreover $D_{+}\subseteq D'_{+},\ \big|W-\Theta\big|_{D'_{+}}\leq 2\varepsilon^{-1} M,$\ hence
$$\big|F_{3}\big|_{{D_{+}}}\leq {{3M}\over {2s}}\big|W-\Theta\big|_{{D_{+}}}\leq{{3M}\over {2s}}\big|W-\Theta\big|_{{D'_{+}}}\leq {{3M}\over {2s}}2\varepsilon^{-1} M\leq q M.$$
Since $q\leq 10^{-2}\theta^{2},\ $ by (I) we get
\begin{equation}\label{e19}
F_{3}\in P^{2}(D_{+}),\ \ \ \big|F_{3}\big|_{{D_{+}}}\leq {1\over 100}\theta^{2} M.
\end{equation}

\medskip
5) Estimate for $F_{1}=w\circ \Omega_{+}-w\circ(\Omega_{+}+\phi)$.
\vspace{0.2cm}

Assume that
\begin{equation}\label{e20}
z\in P^{2}(D_{+}),\ \ \ \big|z\big|_{{D_{+}}}\leq {1\over 24}\theta^{2} M.
\end{equation}
In the previous setting we define
$$\phi={\Theta^{-1}}(z+h^{\star}\circ\Theta).$$
By the definition of $h^{\star}$ we get $\big|h^{\star}\circ \Theta\big|_{D_{+}}\leq M $, and therefore
\begin{equation}\label{e21}
\big|\phi\big|_{D_{+}}\leq \theta^{-1}\Big(\big|z\big|_{D_{+}}+\big|h^{\star}\circ \Theta\big|_{D_{+}}\Big)\leq M\big({\theta\over 24}+\theta^{-1}\big)\\[0.2cm]
 \leq {25\over 24}\theta^{-1}\ M\leq {s\over 100}.
\end{equation}
Thus
\begin{eqnarray*}
\big(\Omega_{+}+\phi\big)\big(D_{+}\big)&\subseteq & D(r_{+}+\varepsilon_{+}s_{+}+{s\over 100},s_{+}+{s\over 100})\\[0.2cm]
&\subseteq& D({r\over 2}+{s\over 2},{11\over 21}s).
\end{eqnarray*}
Since $s<{r\over 9}$, we have
$${r\over 2}+{s\over 2}<{4\over 7}r=r'_{+}+{2\over 3}s-{2\over 21}r<r'_{+}-{1\over 7}s,$$
hence
$$\big(\Omega_{+}+\phi\big)\big(D_{+}\big)\subseteq D({r\over 2}+{s\over 2},{11\over 21}s)\subseteq D^{\star}_{+},$$
$$\big(\Omega_{+}\big)\big(D_{+}\big)\subseteq D(r_{+}+\varepsilon_{+}s_{+},s_{+})\subseteq D({4\over 7}r,{11\over 21}s).$$
From the definition of $r,\ s,\ \theta$, we have
$${s\over \theta}={r\over 300}<{r\over 7},\ \ \ \ {11\over 21}s\leq {2\over 3}s\leq {1\over 3}{s\over \theta},$$
it is easy to see that
$$D({4\over 7}r,{11\over 21}s)\subseteq D(4\rho-{{2s}\over {3\theta}},{s\over {3\theta}})=D(4\rho,t)-{{2s}\over {3\theta}}.$$
Furthermore an application of  Lemma \ref{lem5.1} with $d={{2s}\over {3\theta}}$ yields
$$\big|w(z)-w(z')\big|\leq {{3\theta}\over {2s}}{{2M}\over \varepsilon}|z-z'|\leq q \theta |z-z'|,\ \ \  z,\ z'\in D({4\over 7}r,{11\over 21}s). $$
As a consequence we have
$$\big|F_{1}\big|_{D_{+}}\leq q\theta \big|\phi\big|_{D_{+}}\leq {25\over 24}q M$$
by (\ref{e21}). Since $q\leq 10^{-2} \theta^2$, we have
\begin{equation}\label{e22}
F_{1}\in P^{2}(D_{+}),\ \ \ \big|F_{1}\big|_{D_{+}}\leq  {1\over 96} \theta^2 M.
\end{equation}

\medskip
6) Determination of $z$.
\vspace{0.2cm}

The set of all $z$ satisfying $\big|z\big|_{{D_{+}}}\leq {1\over 24}\theta^{2} M$ is a complete metric space.
Using (\ref{e18}), (\ref{e19}), (\ref{e22}), we get
$$|F(z)|_{D_{+}}\leq \Bigl({1\over 49}+{1\over 100}+{1\over 96}\Bigr)\theta^2 M\leq {1\over 24}\theta^2 M.$$
Hence
$$z \mapsto F(z)=F_{1}+F_{2}+F_{3}\in P^{2}(D_{+})$$
is a mapping of this metric space into itself.

Furthermore letting
$$\phi={\Theta^{-1}}(z+h^{\star}\circ\Theta),\ \phi'={\Theta^{-1}}(z'+h^{\star}\circ\Theta),$$
we have as above
\begin{equation}\label{e30}
\big|F(z)-F(z')\big|_{D_{+}}\leq q\theta \big|\phi-\phi'\big|_{{D_{+}}}=q \theta\big|\Theta^{-1}(z-z')\big|_{D_{+}}\leq q\big|z-z'\big|_{D_{+}}
\end{equation}
for all $z,z'$ satisfying (\ref{e20}), where  $0<q<1$. Hence $F$ is a contraction,\ there is a fixed point $F(z)=z$, which leads to the existence of a mapping $\Phi$ such that the diagram in the inductive theorem exists and commutes. If the unknown function $z$ is represented by $\widetilde{{Z}}(\theta,\eta)$, where $z(\xi,\eta)=\widetilde{Z}(\omega \xi,\eta)$. From the above, we know $\widetilde{{Z}}$ is well defined in $D_+$, and by (\ref{e30}), $\widetilde{{Z}}$ have period $2\pi$ in each of variable $\theta_i (1\leq i\leq n).$ Hence, $z$ is quasi-periodic with the frequency $\omega$ in the first variable which means $\phi$ is quasi-periodic with the frequency $\omega$ in the first variable.

\medskip
7) Proof of inequality (\ref{e8}).
\vspace{0.2cm}

Since $\Phi_{+}-\Omega_{+}=\phi={\Theta^{-1}}(z+h^{\star}\circ\Theta),\ $ we know $\Phi_{+}-\Omega_{+}$ is quasi-periodic with the frequency $\omega$ in the first variable , as a consequence of
$$\big|\Theta^{-1}\circ z\big|_{D_{+}}\leq \theta^{-1}\big|z\big|_{D_{+}},$$
we have the estimate
\begin{eqnarray*}
\big|\Phi_{+}-\Omega_{+}-\Theta^{-1}h^{\star}\circ\Theta\big|_{D_{+}}= \big|\Theta^{-1}\circ z\big|_{D_{+}}\leq \theta^{-1}\big|z\big|_{D_{+}},
\end{eqnarray*}
thus we get
\begin{eqnarray}\label{e23}
\big|\Phi_{+}-\Omega_{+}-\Theta^{-1}h^{\star}\circ\Theta\big|_{D_{+}}\leq \theta^{-1}{1\over 24}\theta^2 M\leq {1\over 24}\theta M.
\end{eqnarray}

The function
\begin{eqnarray}\label{e24}
\eta\mapsto \theta^{-1}[g](\theta\eta)=g_{0}+g_{1}\eta+\cdots,\ g_{j}\in \mathbb{R},\ |\eta|<\theta^{-1}s
\end{eqnarray}
is holomorphic for $|\eta|<\theta^{-1}s$,  and has the estimate $\Big|\theta^{-1}[g](\theta \eta)\Big|\leq \theta^{-1} M,$ since $g$ is the second component of $h$, and $|h(\theta \eta)|\leq M\ \text{for}\ |\eta|<\theta^{-1}s.$

We apply Lemma \ref{lem5.2} to the function (\ref{e24}) with $q={\theta\over 2},\ m=3$ such that for the polynomial
$$Q_{2}(\eta)=a_{0}+a_{1}\eta+a_{2}\eta^2,\ \ \ a_{j}=\Big(1-\big({\theta\over 2}\big)^{6-2j}\Big)g_{j},\ \ \ j=0,1,2,$$
we have the estimate
$$\big|Q_{2}(\eta)-\theta^{-1}[g](\theta\eta)\big|\leq \big({\theta\over 2}\big)^3\theta^{-1}M={1\over 8}\theta^2 M\leq  {1\over 16}\theta M,\ \ \ |\eta|<{\theta\over 2}{s\over \theta}={s\over 2}=s_{+}.$$
Then for
$$Q(\eta)=(0,a_{0}+a_{1}\eta+a_{2}\eta^2),\ \ \ \ \Theta^{-1}\circ h^{\star}\circ\Theta(\eta)=\Big(0,\theta^{-1}[g](\theta\eta)\Big),$$
we obtain $\Phi_{+}-\Omega_{+}-Q$ is quasi-periodic with the frequency $\omega$ in the first variable and
$$
\big|Q-\Theta^{-1}\circ h^{\star}\circ\Theta\big|_{D_{+}}\leq \big|Q_{2}(\eta)-\theta^{-1}[g](\theta\eta)\big|\leq {1\over 16}\theta M.
$$
Together with (\ref{e23}) we have
\begin{eqnarray*}
\big|\Phi_{+}-\Omega_{+}-Q\big|_{D_{+}}&\leq& \big|\Phi_{+}-\Omega_{+}-\Theta^{-1}\circ h^{\star}\circ\Theta\big|_{D_{+}}+\big|Q-\Theta^{-1}\circ h^{\star}\circ\Theta\big|_{D_{+}}\\[0.2cm]
&\leq&  {1\over 24}\theta M+{1\over 16}\theta M={5\over 48}\theta M.
\end{eqnarray*}
The inequality (\ref{e8}) follows. The inductive theorem is proved.\qed

\begin{remark}
Letting
$$\varepsilon=\varepsilon_{k},\ r=r_{k},\ r'=r'_{k},\ s=s_{k},\ s'=s'_{k},\ M=M_{k},$$
$$D=D_{k},\ D'=D'_{k},\ W=W_{k},\ H=H_{k},$$
and replacing  the  index $+$ by $k+1$, Theorem  5.1 confirms what we have asserted  in Section  3 concerning the  construction  of the  commuting diagram (3.12) observing (3.13),\ (3.14),\ (3.15) and (3.16).
\end{remark}

\section{Appendix}
In this appendix we give the detail proof of Lemma \ref{lem2.9} which have been used in the previous sections. For this purpose we need a well known and fundamental approximation result.

\begin{lemma}[Lemma\ 2.1\ in {\cite{Chierchia04}}]\label{lem2.10}
Let $f\in\mathcal{C}^p(\mathbb{R}^\ell)$ for some $p>0$ with finite $\mathcal{C}^p$ norm over $\mathbb{R}^\ell$. Let $\phi$ be a radial-symmetric, $\mathcal{C}^\infty$ function, having as support the closure of the unit ball centered at the origin, also $\phi$ is completely flat and takes
value $1$, let $K=\hat{\phi}$ be its Fourier transform and for all $\delta>0$ define
$$f_{\delta}(x) := K_{\delta}\ast f(x)={\delta}^{-\ell}\int_{\mathbb{R}^\ell} K\Big({{x-y}\over {\delta}}\Big)f(y)dy.$$
Then there exists a constant $c\geq 1$ depending only on $p$ and $\ell$ such that for any $\delta>0$, the function $f_{\delta}(x)$ is real analytic on $\mathbb{C}^\ell$, and
for all $\beta \in \mathbb{N}^{\ell}$ with $|\beta|\leq p$, one has
$$\sup\limits_{x\in \Pi_{\delta}^{\ell}}\Big|\partial^\beta f_{\delta}(x)-\sum\limits_{|\lambda|\leq p-|\beta|} {{\partial^{\lambda+\beta}f(\Re\ x)}\over {\lambda !}}(i\Im\ x)^\lambda \Big|\leq c\big\|f\big\|_{p}{\delta}^{p-|\beta|},$$
and, for all $0<\delta<\delta'$,
$$\sup\limits_{x\in \Pi_{\delta}^{\ell}}\big|\partial^\beta f_{\delta'}-\partial^\beta f_{\delta}\big|\leq c\big\|f\big\|_{p}{\delta'}^{p-|\beta|}.$$
Moreover, the H\"{o}lder norms of $f_{\delta}$ satisfy, for all $0\leq s\leq p\leq r$,
$${\big\|f_{\delta}-f\big\|_{s}}\leq c\big\|f\big\|_{p} {{\delta}^{p-s}},\ \ {\big\|f_{\delta}\big\|_{r}}\leq c \big\|f\big\|_p\delta^{p-r}.$$
The function $f_{\delta}$ preserves periodicity, this is, if $f$ is T-periodic in any of its variable $x_j$, so is $f_{\delta}$.
\end{lemma}

Now we are in a position to prove Lemma \ref{lem2.9}.\\

\noindent\textbf{Proof of Lemma \ref{lem2.9} }\ Since  the function $h(x,y)$ is  quasi-periodic in $x$ with the frequency  $\omega=(\omega_1,\omega_2,\cdots,\omega_n)$,\ from Definition \ref{def2.1}, there exists the corresponding shell function
$$F(\theta,y) :=F(\theta_1,\theta_2,\cdots,\theta_n,y),\ \ \ \theta=(\theta_1,\theta_2,\cdots,\theta_n),$$
which is $2\pi$-periodic in each $\theta_j$, such that
$h(x,y)=F(\omega_1 x,\omega_2 x,\cdots,\omega_n x,y).$

From the assumptions of Lemma \ref{lem2.9}, $h\in \mathcal{C}^p(\mathbb{R}^2)$, then $F\in \mathcal{C}^p(\mathbb{R}^{n+1})$, and $\big\|F\big\|_{p}$ is equivalent to $\big\|h\big\|_{p}$. In fact, if $p\geq 0$ is an integer, then
\begin{align}\label{b8}
\big\|F\big\|_{p}&=\sum\limits_{|\beta|\leq p}\,\sup \limits _{(\theta,y)\in \mathbb{R}^{n+1}}\Big|\partial^\beta F(\theta,y)\Big|=\sum\limits_{|\beta|\leq p}\,\sup \limits _{(\theta,y)\in \mathbb{R}^{n+1}}\Big|{{\partial^\beta F(\theta_1,\cdots,\theta_n,y)}\over{\partial{\theta_{1}^{\beta_{1}}}  \cdots \partial\theta_{n}^{\beta_{n}}  \partial y^{\beta_{n+1}} }}\Big|\nonumber\\[0.2cm]
&=  \sum\limits_{|\beta|\leq p}\,\sup \limits _{(x,y)\in \mathbb{R}^{2}} \Big|{{\partial^\beta h(x,y)}\over {\partial x^{\beta_1+\cdots+\beta_{n}}  \partial y^{\beta_{n+1}} }}\,{1\over{{\omega_1}^{\beta_1}}}\cdots {1\over{{\omega_n}^{\beta_n}}} \Big|\nonumber\\[0.2cm]
 & = \sum\limits_{\lambda_{1}+\lambda_{2}\leq p}\,\sup \limits _{(x,y)\in \mathbb{R}^{2}}\Big|{{\partial^\beta h(x,y)}\over {\partial x^{\lambda_1}  \partial y^{\lambda_{2}} }}\,{1\over{\omega^{\lambda_1}}} \Big|\nonumber\\[0.2cm]
&= \sum\limits_{|\beta|\leq p}\,\sup \limits _{(x,y)\in \mathbb{R}^{2}} {1\over{\big|\omega^{\lambda_1}\big|}}\Big|\partial^\beta h(x,y)\Big| =\tilde{c}\big\|h\big\|_{p}\, .
\end{align}
Similarly,\ if $p=l+s$,\  $l \geq 0$ is an integer,\ $s \in (0,1),\nu:=(x,y),\bar{\nu}:=(\bar{x},\bar{y}),$ then
\begin{eqnarray}\label{b9}
\big\|F\big\|_{p} &=&\sup \limits _{\substack{\nu\neq \bar{\nu}\\ |\beta|=\ell}}{{1\over{|\omega^{\lambda_1}|}}{{\big|\partial^\beta h(\nu)-\partial^\beta h(\bar{\nu})\big|}}\over |\nu-\bar{\nu}|^{s}}+\sum\limits_{|\beta|\leq l}\,\sup \limits _{(x,y)\in \mathbb{R}^{2}} {1\over{\big|\omega^{\lambda_1}\big|}}\Big|\partial^\beta h(x,y)\Big|\nonumber\\[0.2cm]
&= &\tilde{c}\big\|h\big\|_{p},
\end{eqnarray}
where $\beta=(\beta_1,\cdots,\beta_n,\beta_{n+1})\in{\mathbb{N}}^{n+1},|\beta|=\beta_1+\cdots+\beta_{n}+\beta_{n+1},\lambda_1=\beta_1+\cdots+\beta_{n},\lambda_2=\beta_{n+1},\omega^{\lambda_1}={\omega_1}^{\beta_1}\cdots {\omega_n}^{\beta_n},\tilde{c}$ is a positive constant depending only on $p,\omega$.

An application of  Lemma \ref{lem2.10} to the function $F(\theta,y) \in \mathcal{C}^p(\mathbb{R}^{n+1})$, there exists a function $F_{\delta}$, which is an analytic approximation of $F$, and for any $\beta\in{\mathbb{N}}^{n+1}$ with $|\beta|\leq p$,
$$\sup\limits_{z\in \Pi_{\delta}^{n+1}}\Bigg|\partial^\beta F_{\delta}(z)-\sum\limits_{|\lambda|\leq p-|\beta|} {{\partial^{\lambda+\beta}F(\Re\ z)}\over {\lambda !}}(i\Im\ z)^\lambda \Bigg|\leq c_{3}\big\|F\big\|_{p}{\delta}^{p-|\beta|}$$
and for any $0<\delta<\delta',$
$$\sup\limits_{z\in \Pi_{\delta}^{n+1}}\big|\partial^\beta F_{\delta'}-\partial^\beta F_{\delta}\big|\leq c_{3}\big\|F\big\|_{p}{\delta'}^{p-|\beta|},$$
where $c_{3}\geq 1$ is a positive constant depending only on $p$ and $n$.

Specially, if $\beta=0$, one has
$$\sup\limits_{z\in \Pi_{\delta}^{n+1}}\Bigg|F_{\delta}(z)-\sum\limits_{|\lambda|\leq p} {{\partial^{\lambda}F(\Re\ z)}\over {\lambda !}}(i\Im\ z)^\lambda \Bigg|\leq c_{3}\big\|F\big\|_{p}{\delta}^{p},$$
$$\sup\limits_{z\in \Pi_{\delta}^{n+1}}\big|F_{\delta'}-F_{\delta}\big|\leq c_{3}\big\|F\big\|_{p}{\delta'}^{p},$$
where $0<\delta<\delta'$.
Hence,
\begin{eqnarray*}
\sup\limits_{(\theta,y) \in \Pi_{\delta}^{n+1}}\big|F_{\delta}(\theta,y)\big| &\leq & \sup\limits_{z\in \Pi_{\delta}^{n+1}}\Bigg|\sum\limits_{|\lambda|\leq p} {{\partial^{\lambda}F(\Re\ z)}\over {\lambda !}}(i\Im\ z)^\lambda \Bigg|+ c_{3}\big\|F\big\|_{p}{\delta}^{p} \\[0.2cm]
 & \leq& c_{4}\big\|F\big\|_{p}+ c_{3}\big\|F\big\|_{p}{\delta}^{p},
\end{eqnarray*}
where $c_{4}$ is a positive constant depending only on $p,n.$

Now we define the analytic approximation of $h(x,y)$ in $E_\delta$ as follows
$$h_{\delta}(x,y)=F_{\delta}(\omega_1 x,\omega_2 x,\cdots,\omega_n x, y).$$
Thus
$$\big|h_{\delta}\big|_{E_\delta} :=\sup\limits_{(\theta,y) \in \Pi_{\delta}^{n+1}}\big|F_{\delta}(\theta,y)\big|\leq c_{4}\big\|F\big\|_{p}+ c_{3}\big\|F\big\|_{p}{\delta}^{p}.$$
By (\ref{b8}), (\ref{b9}), we have
$$\big|h_{\delta}\big|_{E_\delta}\leq c_{4}\big\|F\big\|_{p}+ c_{3}\big\|F\big\|_{p}{\delta}^{p}\leq  c_{2}\big\|h\big\|_{p}.$$
Similarly,  for any $0<\delta<\delta',$ one can obtain
$$\big|h_{\delta'}-h_{\delta}\big|_{E_\delta}\leq \sup\limits_{z\in \Pi_{\delta}^{n+1}}\big|F_{\delta'}-F_{\delta}\big|\leq c_{2}\big\|h\big\|_{p}{\delta'}^{p},$$
where $c_{2}>0$ is a constant depending only on $p,n,\omega.$ Hence,
$$\big|h_{\delta'}-h_{\delta}\big|_{E_\delta}\leq c_{2}\big\|h\big\|_p{\delta'}^{p}.$$
The proof of this lemma is completed. \qed

\bigskip

\section*{References}
\bibliographystyle{elsarticle-num}

\begin{thebibliography}{99}
{\small
\rm


\bibitem{Chierchia04} L. Chierchia, D. Qian, {\em Moser's theorem for lower dimensional tori},
J. Differential Equations 206 (2004) pp. 55-93.

\bibitem{Herman83} M.\ R.\ Herman, {\em  Sur les courbes invariantes par les diff\'{e}omorphismes de l'anneau I}, Ast\'{e}risque No. 103-104 (1983).

\bibitem{Herman86} M.\ R.\ Herman, {\em  Sur les courbes invariantes par les diff\'{e}omorphismes de l'anneau II}, Ast\'{e}risque No. 144 (1986).

\bibitem{Huang16} P. Huang, X. Li, B. Liu {\em Quasi-periodic solutions for an asymmetric oscillation}, Nonlinearity 29 (2016) pp. 3006-3030.

\bibitem{Levi01} M.\ Levi, J.\ Moser, {\em A Lagrangian proof of the invariant curve theorem for twist mappings Smooth Ergodic Theory and its Applications}, (Seattle, WA, 1999) (Proc. Symp. Pure Math. vol 69) (Providence, RI: American Mathematical Society) (2001) pp. 733-46.

\bibitem{Liu05} B.\ Liu, {\em Invariant curves of quasi-periodic reversible mapping}, Nonlinearity 18 (2005) pp. 685-701.

\bibitem{Moser62} J.\ Moser,  {\em On invariant curves of area-perserving mappings of an annulus}, Nachr. Akad. Wiss. G\"{o}ttingen Math. -Phys.  vol II (1962) pp. 1-20.

\bibitem{Moser66} J.\ Moser,  {\em  A Rapidly Convergent Iteration  Method and Nonlinear Differential
Equations II}, Ann.Scuola  Norm. Sup. Pisa (1966) pp. 499-535.

\bibitem{Moser88} J.\ Moser, {\em A stability theorem for minimal foliations on a torus Ergod}, Theory Dynam. Syst.\ 8 (1988) pp. 251-81.

\bibitem{Russmann70} H.\ R\"{u}ssmann, {\em  Kleine Nenner I: \"{U}ber invariante Kurven differenzierbarer Abbildungen eines Kreisringes}, Nachr. Akad. Wiss. G\"{o}ttingen Math.-Phys. Kl. II (1970) pp. 67-105.

\bibitem{Russmann74} H.\ R\"{u}ssmann, {\em  On optimal estimates for the solutions of linear partial differential equations of first order with constant coefficients on the torus}, Dynamical systems theory and applications (1974) pp. 598-624.

\bibitem{Russmann83} H.\ R\"{u}ssmann, {\em  On the existence of invariant curves of twist mappings of an annulus},  Lecture Notes in Math.  Springer Berlin 1007 (1983) pp. 677-718.

\bibitem{Siegel97} C.\ Siegel  and J.\ Moser, {\em Lectures on Celestial Mechanics}, (Berlin: Springer) (1997).

\bibitem{Zharnitsky00} V.\ Zharnitsky, {\em  Invariant curve theorem for quasiperiodic twist mappings and stability of motion in the Fermi-Ulam problem}, Nonlinearity 13 (2000) pp. 1123-36.

\bibitem{Zehnder75} E.\ Zehnder, {\em  Generalized Implicit  Function Theorems with Applications to  Some
Small Divisor Problems}, I.  Comm. Pure Appl. Math.  28 (1975) pp. 91-140.


}

\end{thebibliography}

\end{document}